\documentclass[12pt]{article}
\usepackage{amsmath,latexsym,epsfig}
\usepackage{amscd}
\usepackage{mathrsfs,amsfonts,stmaryrd,bbm}
\usepackage{amssymb}

\DeclareMathSymbol{\bdC}{\mathbin}{AMSb}{'103}
\DeclareMathAlphabet{\mathpzc}{OT1}{pzc}{m}{it}
\DeclareMathSymbol{\subsetne}{\mathbin}{AMSb}{'050}

\newtheorem{thm}{Theorem}[section]
\newtheorem{lem}[thm]{Lemma}
\newtheorem{cor}[thm]{Corollary}
\newtheorem{prob}[thm]{Problem}
\newtheorem{defn}[thm]{Definition}

\newtheorem{example}[thm]{Example}

\numberwithin{equation}{section}

\newcommand{\cP}{\mathscr{P}}
\newcommand{\mP}{\mathnormal{P}}
\newcommand{\cPjmn}{\cP^j_{m,n}}
\newcommand{\cPkmn}{\cP^k_{m,n}}
\newcommand{\cPmn}{\cP_{m,n}}
\newcommand{\mPmn}{\mP_{m,n}}

\newcommand{\C}{\bdC}
\newcommand{\h}{{{\mbox{\tiny $\mathsf{H}$}}}}
\newcommand{\spn}{\mathpzc{span}}
\newcommand{\dist}[1]{\mathpzc{dist}\left(\,#1\,\right)}
\newcommand{\mT}{\mathbb{T}}
\newcommand{\mS}{\mathbb{S}}
\newcommand{\nullity}[1]{\mathpzc{nullity}\left(\,#1\,\right)}

\newcommand{\sg}{\sigma}
\newcommand{\sgn}[1]{\sigma_{-#1}}
\newcommand{\sgmin}{\sgn{1}}
\newcommand{\dm}{\mathpzc{dim}}
\newcommand{\dg}{\mathpzc{deg}}
\newcommand{\al}{\alpha}
\newcommand{\bt}{\beta}

\newcommand{\dl}{\delta}
\newcommand{\GCD}{\mathpzc{gcd}}
\newcommand{\prf}{\noindent {\bf Proof. \ }}
\newcommand{\cL}{\mathcal{L}}
\newcommand{\codim}{\mathpzc{codim}}
\newcommand{\eps}{\varepsilon}
\newcommand{\cD}{\mathcal{D}}
\newcommand{\Dl}{\Delta}

\newcommand{\bdu}{\mathbf{u}}
\newcommand{\bdv}{\mathbf{v}}
\newcommand{\bdw}{\mathbf{w}}
\newcommand{\bdx}{\mathbf{x}}
\newcommand{\bdy}{\mathbf{y}}
\newcommand{\bdz}{\mathbf{z}}
\newcommand{\bda}{\mathbf{a}}
\newcommand{\bdb}{\mathbf{b}}
\newcommand{\bdc}{\mathbf{c}}

\newcommand{\bdf}{\mathbf{f}}
\newcommand{\bdg}{\mathbf{g}}
\newcommand{\bdh}{\mathbf{h}}

\newcommand{\bdk}{\mathbf{k}}

\newcommand{\bdp}{\mathbf{p}}
\newcommand{\bdq}{\mathbf{q}}

\newcommand{\bdo}{\mathbf{0}}

\newcommand{\qed}{\hfill $\blacksquare$}

\begin{document}

\title{The Numerical Greatest Common Divisor
of Univariate Polynomials}

\author{Zhonggang Zeng \thanks{
Department of Mathematics, Northeastern Illinois University,
Chicago, IL 60625, Email: {\tt zzeng@neiu.edu}, 
Research supported in part by NSF under Grant DMS-0715127}
}


\date{July 4, 2011}


\maketitle

\begin{abstract}
This paper presents a regularization theory for numerical
computation of polynomial greatest common divisors and a
convergence analysis, along with a detailed description of
a blackbox-type algorithm.
~The root of the ill-posedness in conventional GCD computation is 
identified by its geometry where polynomials form differentiable 
manifolds entangled in a stratification structure.
~With a proper regularization, 
the numerical GCD is proved to be strongly well-posed.
~Most importantly, the numerical GCD solves the problem of finding the
GCD accurately using floating point arithmetic even if the data are
perturbed.
~A sensitivity measurement, error bounds at each computing stage, and 
the overall convergence are established rigorously.
~The computing results of selected test examples show that the algorithm 
and software appear to be robust and accurate.
\end{abstract}

\section{Introduction}

As one of the fundamental algebraic problems with a long history,
finding the greatest common divisor (GCD) of univariate polynomials 
is an indispensable component of many algebraic computations besides
being an important problem in its own right.
~The classical Euclidean Algorithm has been known for centuries
\cite[p.58]{GatGer} and the problem is well studied in computer algebra,
~where algorithms are developed using exact arithmetic with exact
data.
~These algorithms are not suitable for practical numerical computation
because computing GCD is an ill-posed problem in the sense that it is
infinitely sensitive to round-off error and data perturbations.
~A tiny error in coefficients generically degrades the GCD into a
meaningless constant.
~The central problem of this paper is: ~How, and why, can we still recover
the lost GCD accurately using the inexact data and floating point 
arithmetic?

For this purpose, we study the root of the ill-posedness by presenting
the geometry of the polynomial GCD problem: ~The collection of polynomial
pairs whose GCD's share a fixed positive degree forms a differential 
manifold of a positive codimension, and these manifolds entangle in
a stratification structure in which a manifold is in the closure of
manifolds of lower codimensions.
~The ill-posedness of the GCD problem lies in the dimension deficit
of the GCD manifold from which a polynomial pair will pushed away by
arbitrary perturbations.

Taking advantage of the geometric properties, we study the 
{\em numerical GCD} formulated by Corless, Gianni, Trager and Watt
\cite{corless-gianni} as well as Karmarkar and Lakshman \cite{kar-Lak96}
by establishing a comprehensive regularization theory of numerical GCD.
~We prove that numerical GCD generalizes the traditional notion of GCD
and, when the data are sufficiently accurate, the numerical GCD uniquely
exists and is Lipschitz continuous, thereby making it strongly well-posed 
and computable using floating point arithmetic.
~Most importantly, the numerical GCD solves the central problem
of this paper by approximating the exact GCD of the underlying
polynomials hidden in data perturbation.

Building upon the thorough study on the theory of numerical GCD, we
further establish a detailed analysis of the algorithm proposed by
the author \cite{zeng-mr-05} for computing the numerical GCD,
and prove the Numerical GCD Convergence Theorem.
~The algorithm is implemented using the code name {\sc uvGCD}
as part of the comprehensive package {\sc ApaTools} \cite{apatools}
for approximate polynomial algebra.

As mentioned above, GCD-finding is one of the basic operations in
algebraic computation with a wide range of applications that include
engineering problems such as graphics and modeling, robotics, computer
vision, image restoration, control theory, system identification
\cite{barnett,fau93,hen-seb,kar-mit-00,mer90,lia-pil,sed-chang,sto-so,
ver-wang},
as well as other branches of mathematics and computer science such as
simplifying rational expressions, partial fraction expansions, 
canonical transformations, mechanical geometry theorem proving
\cite{chousc,geddes,zippel}, hybrid rational function approximation
\cite{kai-nod}, and decoder implementation for error-correction
\cite{bre-kun}.
~A robust GCD-finder is also crucial to
root-finding studies when multiple roots are present
\cite{dun,noda-sas,zeng-mr-05}.
~In recent years, substantial effort has been spent on developing
algorithms for computing the numerical GCD of inexact polynomials.
~These pioneering works include resultant-based algorithms
\cite{corless-gianni,emiris-galligo-lombardi,rupprecht},
optimization strategies
\cite{chin-corless,kar-Lak96},
modifications of the Euclidean Algorithm
\cite{beck-fast,hribernig-stetter,noda-sas,schonhage},
root grouping \cite{pan96,ver-wang}, QR factorization
\cite{cor-watt-zhi,zaro}, and low rank approximations
\cite{kmyz05,KalYanZhi06,WinAll08b,WinAll08a}.
~Several methods have been implemented as part of Maple
SNAP package \cite{jean-lab} that include {\sc QuasiGCD}, {\sc EpsilonGCD}
and {\sc QRGCD}.
~Particularly in \cite{corless-gianni}, Corless, Gianni, Trager and Watt
propose a novel, albeit unfinished, approach that includes the use of
singular value decomposition to identify the GCD degree and
several possible strategies for calculating the GCD factors
including solving least squares problems.

In the context of polynomial root-finding, we developed
a new special case algorithm for computing the GCD of a polynomial and its
derivative.
~This algorithm is briefly described in \cite{zeng-mr-05}
as an integral component of Algorithm {\sc MultRoot} \cite{zeng_multroot} that
calculates multiple roots of a polynomial with high accuracy without
using multiprecision arithmetic even if the coefficients are perturbed.
~Due to the scope of the paper \cite{zeng-mr-05}, that algorithm is narrowly
featured without in-depth analysis of the problem regularity, algorithm
convergence, error analysis, extensions, applications, or comprehensive
testing/experiment.
~Our numerical GCD algorithm employs a successive Sylvester matrix updating
process
for identifying the maximum degree of the numerical GCD along with
an initial approximation to the GCD factors.
~Then the Gauss-Newton iteration is applied to certify the GCD and to
refine the polynomial factors via solving a regular quadratic least
squares problem.
~Those new strategies apparently fill the main gaps in
previous works and is mentioned in a recently published textbook:

\begin{quote}
``This numerical common gcd
algorithm ... appears to be the most efficient and reliable algorithm
for that purpose; I have seen it too late to include it in the text.''
~~H. J. Stetter, {\em Numerical Polynomial Algebra} \cite[p.223]{stetter_book}
\end{quote}

The software {\sc uvGCD} has been tested rigorously and extensively. As
sample test results shown in \S \ref{sec:res},
{\sc uvGCD} is substantially more robust and accurate than the existing
packages.
~The complexity of our method is ~$O(n^3)$~ for
the combined degree ~$n$~ of the given polynomials.

The main theorems in this paper appear to be new, including GCD Manifold
Theorem, GCD Extension Theorem, Numerical GCD Regularity Theorem, 
Numerical GCD Approximation Theorem and Numerical GCD Convergence 
Theorem.

\section{Difficulties of finding GCD in numerical computation}

Computing polynomial GCD is a typical ``ill-posed problem'' whose 
numerical solutions are generally unattainable using conventional methods, 
even if the method is among the most celebrated in the history.
~The hypersensitivity of such problems can be illustrated in a simple
example:

\begin{example} \label{ilex} Consider a pair of polynomials
\begin{eqnarray*}
   p(x)  &=& \mbox{$x^{10}+ \frac{31}{3} x^9+\frac{10}{3}x^8+x+10$} \\
   q(x)  &=& \mbox{$x^{10}+\frac{71}{7}x^9+\frac{10}{7}x^8 -\frac{6}{7}x
-\frac{60}{7}$}
\end{eqnarray*}
They can be factored as ~$(x+10)(x^9+\frac{x^8}{3}+1)$ ~and 
~$(x+10)(x^{9}+\frac{x^8}{7}-\frac{6}{7})$ ~respectively.
~There is no difficulty for a common computer algebra system like Maple
to find the GCD using symbolic computation:
\end{example}

\footnotesize
\begin{verbatim}
     > gcd(x^10+(31/3)*x^9+(10/3)*x^8+x+10,x^10+71*x^9/7+10*x^8/7-6*x/7-60/7);
\end{verbatim}
\[ x+10\]
\normalsize

However, the GCD quickly degrades from ~$x+10$ ~to a constant simply by
replacing the fractional numbers in the coefficients with
floating point values at the simulated hardware precision of 10 digits%
\footnote{The test is carried out using Maple 12. ~Other versions of
Maple may yield different results}:

\footnotesize
\begin{verbatim}
     > gcd(x^10+10.33333333*x^9+3.333333333*x^8+x+10.,
           x^10+10.14285714*x^9+1.428571429*x^8-.8571428571*x-8.571428571);
\end{verbatim}
\[ 1.000000000\]
\normalsize

The constant 1 is, as a matter of fact, the correct GCD in exact sense 
from the given
(perturbed) coefficients, and the nontrivial GCD ~$x+10$ ~is lost from a
tiny perturbation in data.
~This is the ill-posed nature of GCD for being discontinuous with respect to
its coefficients.

For hundreds of years, the classical Euclidean Algorithm has been {\em the}
method for GCD finding.
~However, it can easily fail in numerical computation.
~The Euclidean Algorithm is a recursive process of polynomial division
\begin{equation} \label{pdiv}
  f ~=~ q\cdot g + r
\end{equation}
from a polynomial pair ~$(f,g)$ ~to obtain the quotient ~$q$ ~and the
remainder ~$r$.
~Assume the degree of ~$g$ ~is no larger than that of ~$f$ ~without loss of
generality, and initialize ~$f_0 = f$, ~$f_1 = g$.
~The Euclidean Algorithm
\begin{equation} \label{eucalg}
 f_j ~=~ q_j \cdot f_{j+1} + f_{j+2},  ~~~~j = 0, 1, \ldots
\end{equation}
generates a remainder sequence ~$f_2, ~f_3, ~\ldots$ ~that terminates at
a scalar multiple of the GCD of ~$f$ ~and ~$g$.
~The Euclidean Algorithm and its modifications remain the method of choice
for GCD computation in symbolic computation for exact polynomials.
~The following example illustrates why the Euclidean Algorithm
behaves poorly in the presence of data error or round-off.

\begin{example} \label{euex} ~~Consider polynomials appeared in
Example~\ref{ilex};
\begin{equation} \label{fg}
\begin{array}{rcl}
   f(x)  &=& \mbox{$(x+10)(x^9+\frac{x^8}{3}+1)$} ~=~ 
\mbox{$x^{10}+ \frac{31}{3} x^9+\frac{10}{3}x^8+x+10$},  \\
   g(x)  &=& x + 10.
\end{array}
\end{equation}
with their GCD equals to ~$g(x)$.
~The (exact) Euclidean Algorithm stops at one step since
~$f(x) = (x^9+\frac{x^8}{3}+1)\cdot g(x) + 0$.
~However, if the coefficients of ~$f$ ~are inexact with a perturbation of
a magnitude at the simulated hardware precision,
say 
\[ \tilde{f}(x) = x^{10}+10.33333333 x^9+3.333333333 x^8+x+10.
\]
the one-step Euclidean Algorithm involves a polynomial division
\begin{eqnarray*} 
\tilde{f}(x) &=& \mbox{\scriptsize $(x+10)(
x^9+.33333333 x^8+.000000033 x^7-.00000033 x^6+
$} \\
& & \mbox{\scriptsize $ 
+.0000033 x^5-.000033 x^4+.00033 x^3-.0033 x^2+.033 x+0.67) +3.3$}
\end{eqnarray*}
Far from getting a zero, the remainder becomes $3.3$, and
the Euclidean Algorithm produces a constant GCD with a large difference
from ~$x+10$ ~we are looking for,  even though
the data error is tiny. \qed
\end{example}

The numerical instability of the Euclidean Algorithm is inherent from
polynomial division (\ref{pdiv}), which is equivalent to linear system
for the coefficients of $q$ ~and $r$.
Using the $f$ ~and $g$ in (\ref{fg}) as an example again, equation
(\ref{pdiv}) can be written as
\begin{equation} \label{pdsys}
 \mbox{\scriptsize $\left[ \begin{array}{rrrrr} 1 & & & &   \\
10 & 1 & & & \\ & 10 & \ddots & & \\ & & \ddots & 1 &  \\
& & & 10 & 1 \end{array} \right] $}
\mbox{\footnotesize $
\left[ \begin{array}{c} q_{9} \\ q_{8} \\ \vdots \\ q_0 \\ r_0 \end{array}
\right]
$} ~~=~~
\mbox{\scriptsize $
\left[ \begin{array}{c} 1 \\ 10.33333333 \\ 3.333333333 \\ 0 \\ \vdots \\ 0 
\\ 1 \\ 10 \end{array} \right]
$}
\end{equation}
Perturbations in coefficients of ~$f$ ~and/or ~$g$ ~are magnified by
the large condition number ~$1.1 \times 10^{10}$ ~of the ~$10\times 10$
~matrix in (\ref{pdsys}), as shown by the nonzero remainder computed
in Example~\ref{euex}.

The question is: Can we accurately compute the GCD, say $x+10$ ~of
$p(x)$ and $q(x)$ in Example~\ref{ilex}, using the inexact data and
floating point arithmetic?
~The examples above indicate the futility of computing the {\em exact}
GCD in this situation.
~Instead, we compute the {\em numerical greatest common divisor} to
be formulated in \S \ref{sec:agcd}.
~Such a numerical GCD will be proven insensitive with a finite GCD condition
number.
~More importantly, the numerical GCD approximates the underlying exact
GCD, say ~$x+10$, ~with an error in the order of data perturbation,
as confirmed by our Maple software {\sc uvGCD} result:

\footnotesize
\begin{verbatim}
    > u, v, w, res := uvGCD(x^10+10.33333333*x^9+3.333333333*x^8+x+10, 
        x^10+10.14285714*x^9+1.428571429*x^8-.8571428571*x-8.571428571,x,1E-8):
    > u/lcoeff(u,x)   # scale the numerical GCD to be monic
\end{verbatim}
\[9.999999998+1.000000000\,x  \]
\normalsize
The numerical GCD within $10^{-8}$ is
$\tilde{u} = x + 9.999999998$, an accurate approximation to the
exact GCD $x + 10$.

\section{Preliminaries} \label{sec:prelim}

The fields of complex numbers are denoted by
~$\C$.\label{'RC'}
~The ~$n$~ dimensional complex vector space is denoted by ~$\C^n$, ~in which
vectors are columns denoted by boldface
lower case letters such as ~$\bda$, ~$\bdu$, ~$\bdv_2$, ~etc, with ~$\bdo$
~being a zero vector whose dimension can be understood from the context.
\label{'vector'}
~Matrices are represented by upper case letters like ~$A$ ~and ~$J$.
\label{'matrix'}
~For every vector or matrix ~$(\cdot)$, ~the notation
~$(\cdot)^\top$~ represents the transpose and
~$(\cdot)^\h$~ the Hermitian adjoint (i.e. conjugate transpose) of ~$(\cdot)$.
\label{'top'}
~We find it convenient to use Matlab notation ``;'' to stack (column)
vectors as
\[
[3;\,-2;\,4] ~\equiv~
\mbox{\scriptsize $
\left[ \begin{array}{r} 3 \\ -2 \\ 4 \end{array} \right]$}
~\equiv~ [3,\,-2,\,4]^\top, ~~~~
[\bdu;\bdv] \;\equiv\; \mbox{\small
$\left[\begin{array}{c} \bdu \\ \bdv \end{array} \right]$} \;\equiv\;
[\bdu^\top,\bdv^\top]^\top,
\]
The norm $\|\bdv\|$ of a vector $\bdv$ is the Euclidean norm
$\|\bdv\| = \sqrt{\bdv^\h\bdv}$ ~throughout this paper.\label{'vnorm'}
The matrix norm $\|A\|$ of $A$ is induced from the vector
norm $\displaystyle \|A\| \;=\; \max_{\|\bdx\|=1} \|A\bdx\|$. \label{'mnorm'}
We also use the Frobenius norm \cite[Page 55]{golub-vanloan} denoted by
$\|\cdot\|_F$ in some occasions.

All vector spaces are in ~$\C$.
~A vector space spanned by vectors ~$\bdv_1,\ldots,\bdv_n$ ~is denoted by
~$\spn\{\bdv_1,\ldots,\bdv_n\}$.
~The notation ~$\dist{\mS,\mT}$ ~stands for the distance between two
subspaces ~$\mS$ ~and ~$\mT$ ~in a larger vector space
\cite[p. 76]{golub-vanloan}.
~The dimension of the kernel of matrix ~$A$ ~is
~$\nullity{A}$, ~namely the nullity of ~$A$.

For any matrix ~$A$~ of ~$m \times n$~ with
~$m\geq n$, ~there are ~$n$~ singular values \cite[\S 2.5.3]{golub-vanloan}
\[ \sg_1(A) ~\ge~ \sg_2(A) ~\ge~ \ldots ~\ge~ \sg_n(A) ~\ge~ 0. \]
of ~$A$~ with ~$\sg_1(A) = \|A\|$. \label{'sgj'}
We shall also denote the same singular values in reversed order \label{'sgnj'}
\[ 0 ~\le~ \sgn{1}(A) ~\le~ \sgn{2}(A) ~\le~ \cdots ~\le~ \sgn{n}(A) ~=~
\|A\|. \]
Singular value ~$\sgn{1}(A)$~ is the smallest distance from ~$A$~ to a
matrix that is rank-deficient by one. \label{'sgm'}
~Likewise, singular value ~$\sgn{2}(A)$~ is the smallest distance from
~$A$~ to a matrix that is rank-deficient by two, and so on.
~The matrix 
\[ A^+ = (A^\h A)^{-1} A^\h \]
exists uniquely as the
Moore-Penrose inverse of ~$A$~ when ~$\sgn{1}(A) > 0$.  \label{'+'}
~It is straightforward to verify that
\[
\sgmin(A) \;=\; \left\| A^+ \right\|^{-1} \;=\; \min_{\|\bdx\|=1} \|A\bdx\|.
\]
This minimum is attainable at the right singular vector
~$\bdy$~ of ~$A$~ corresponding to ~$\sgmin(A)$.
~Namely ~$\left\|A\bdy\right\| \,=\,\sgmin(A)$.

In this paper, polynomials are in ~$\C$ ~in a single variable
~$x$.
~The ring of such polynomials is denoted by ~$\C[x]$.
~A polynomial is denoted by lower a case
letter, say ~$f$, ~$v$, ~or ~$p_1$, ~etc.
~A polynomial
\begin{equation} \label{pp}
 p ~=~ \rho_0 + \rho_1 x + \rho_2 x^2 + \cdots + \rho_n x^n
\end{equation}
is of degree ~$n$ ~if ~$\rho_n \ne 0$, ~or the degree is ~$-\infty$
~if ~$f(x) \equiv 0$.
~We shall denote the degree of a polynomial ~$p$ ~by ~$\dg(p)$.
\label{'deg'}

For an integer ~$n$, ~the collection
of polynomials with degrees less than or equal to ~$n$
~form a vector space  \label{'mpbdn'}
\[ \mP_n ~~=~~
\big\{ p \in \C[x] ~\big|~ \dg(p) \le n \big\}.
\]
Thus the dimension of ~$\mP_n$ ~is\label{'dim'} \label{'nnbdn'}
\[
\dm(\mP_n) ~~=~~ \left\{ \begin{array}{ccl} 0 & & \mbox{for}~~ n < 0 \\
n + 1 & & \mbox{for}~~ n \ge 0. \end{array} \right.
\]

Throughout this paper, we use the monomial basis ~$\{1,x,x^2,\ldots,x^n\}$
~for ~$\mP_n$, ~in which
every polynomial ~$p$ ~can be written in the form of (\ref{pp})
and corresponds to a
{\em coefficient vector}
\[ \bdp ~=~ [\rho_0;~\ldots;~\rho_n] ~\in~ \C^{n+1}. \]
Notice that ~$\rho_n = 0$ ~is possible, and the same polynomial ~$p$ ~can
be embedded in the space ~$\mP_m \supset \mP_n$ ~with a coefficient vector 
of higher dimension.
~Throughout this paper, if a letter (say ~$f$, ~$g$, ~$q_1$) represents a
polynomial, the same letter in boldface, say ~$\bdf$,
~$\bdg$, ~$\bdq_1$, ~is its coefficient vector
within a vector space ~$\mP_n$ ~that is clear from the context.
~The norm ~$\|a\|$~ of polynomial ~$a \in \mP_n$~ is
defined as the Euclidean norm ~$\|\bda\|$ ~of its coefficient vector ~$\bda$.

We denote the vector space of polynomial pairs
~$(p,q) \in \mP_m \times \mP_n$ ~as
\[ \mPmn  ~=~
\big\{ (p,q) \in \C[x]^2 ~\big|~ \dg(p) \le m, ~\dg(q) \le n \big\},
\]
and its subset formed by polynomial pairs of degrees equal to ~$m$
~and ~$n$, ~respectively, as
\[
\cPmn ~=~
\big\{ (p,q) \in \mPmn ~\big|~ \dg(p) = m, ~\dg(q) = n \big\}.
\]

For every polynomial pair ~$(p,q)$, ~a {\em greatest common divisor} or
GCD of ~$(p,q)$ ~is {\em any} polynomial ~$u$ ~of the highest degree that 
divides both ~$p$ ~and ~$q$.
~Notice that we do not require a GCD to be monic here to avoid
scaling a polynomial by a tiny leading coefficient in computation.
~In this setting, GCD's are not unique and two GCD's of the same polynomial
pair differ by a nonzero constant multiple.
~We define an equivalence relation ~$\sim$ ~between two polynomials in the
sense that ~$f \sim g$ ~if ~$f = \al\, g$ ~for a constant ~$\al \ne 0$.
~Thus the collection of all GCD's of a polynomial pair ~$(p,q)$ ~forms
a ~$\sim$-equivalence
class, denoted by ~$\GCD(p,q)$, ~which is unique in the quotient ring
~$\C[x]/\sim$.

The collection of polynomial pairs with a specified GCD degree ~$k$ ~is
denoted by
\begin{equation} \label{cpbdkmn}
\cPkmn ~=~
\big\{ (p,q) \in \cPmn ~\big|~ \dg\big(\GCD(p,q)\big) = k \big\}
\end{equation}
If ~$u \in \GCD(p,q)$, then polynomials ~$v = p/u$ ~and
~$w = q/u$ ~are called the {\em cofactors} of polynomial pair ~$(p,q)$.
~The distance between two polynomial pairs, or generally the distance
between two polynomial arrays ~$(p_1,\ldots,p_l)$~ and
~$(q_1,\ldots,q_l)$~ is naturally derived from the polynomial norm
\begin{equation} \label{arraynorm}
 \big\|(p_1,\ldots,p_l)-(q_1,\ldots,q_l)\big\| ~~=~~
\sqrt{\|p_1-q_1\|^2 + \ldots + \|p_l - q_l\|^2}.
\end{equation}

Let $\Psi_n ~:~ \mP_n \longrightarrow \C^{n+1}$ denote the
isomorphism that maps a polynomial $a$ in $\mP_n$ ~to its coefficient
vector $\bda$ in $\C^{n+1}$, namely $\Psi_n(a) = \bda$. \label{'psin'}
~For a fixed ~$f \in \mP_n $ ~and any ~$g \in \mP_m$,
~the polynomial multiplication ~$f\,g$ ~is a linear transformation
\[ \mathcal{F}_m \,:\, \mP_m  \longrightarrow \mP_{m+n}
\mbox{\ \ \ with \ \ \ } \mathcal{F}_m(g) = f\cdot g
~~~\mbox{for every} ~~g \in \mP_m.
\]
Let ~$\bdf = [\phi_0;~\phi_1;~\ldots;~\phi_n]$ ~be the coefficient vector
of ~$f$.
~The matrix for the linear transformation ~$\mathcal{F}_m$ ~is called a
convolution matrix
\begin{equation} \label{cauchy}
 C_m(f) \;=\; \mbox{\scriptsize $
\begin{array}{c}
\overbrace{\mbox{\ \ \ \ \ \ \ \ \ \ \ \ \ \ \ \ \ \ }}^{m+1}
\\
\left[
\begin{array}{ccc}
\phi_0 & &  \\
\vdots & \ddots & \\
\phi_n &        & \phi_0 \\
&  \ddots  & \vdots \\
& & \phi_n
\end{array} \right]
\end{array}$}.
\end{equation}
For polynomials ~$u \in \mP_j$~ and ~$v \in \mP_k$~ with coefficient
vectors ~$\bdu$ ~and ~$\bdv$ ~respectively,
\[ \bdw ~=~ C_k(u) \,\bdv ~=~ C_j(v) \,\bdu ~=~ \Psi_{j+k}(u\cdot v) \]
is the coefficient vector ~$\bdw \in \C^{j+k+1}$ ~of polynomial
product ~$w = u \, v \in \mP_{j+k}$.

The classical Sylvester matrices can be derived naturally.
~Let ~$(p,q)$ ~be a given pair of polynomials of degrees ~$m$ ~and ~$n$
~respectively, if ~$u$ ~is a GCD of ~$(p,q)$ ~with cofactors ~$v$ ~and ~$w$.
~Then ~$p \cdot w - q \cdot v \,=\, uvw-uwv \,=\, 0$, ~namely
\begin{equation} \label{bdpwqv0}
   C_{n-j}(p) \, \bdw - C_{m-j}(q) \, \bdv ~~\equiv~~
   \big[ C_{n-j}(p) ~|~ C_{m-j}(q) \big] \,
\left[ \begin{array}{r} \bdw \\ -\bdv \end{array}
\right] ~~=~~ \bdo.
\end{equation}
for any degree ~$j \le \dg(u)$.
~In other words, matrix
~$\big[ C_{n-j}(p) ~|~ C_{m-j}(q) \big]$ ~is rank-deficient
if ~$j \le \dg(u)$, ~and the GCD problem is equivalent to the rank/kernel
problem of such matrices.

With $p(x) = p_0 + p_1 x + \ldots + p_m x^m \in \mP_m$,
$q(x) = q_0 + q_1 x + \cdots + q_n x^n \in \mP_n$, and $j = 1, 2, \ldots,
\min\{m,n\}$,
the $j$-th {\em Sylvester matrix} \,of $(p,q)$  in $\mPmn$ ~is
defined as
\begin{equation} \label{Sj}
 S_j(p,q) ~~=~~ \big[ C_{n-j}(p) ~|~ C_{m-j}(q) \big]
~~=~~
\begin{array}{cc} \;\;\;\;
\overbrace{\mbox{\ \ \ \ \ \ \ \ \ \ \ \ \ \ }}^{n-j+1} &
\overbrace{\mbox{\ \ \ \ \ \ \ \ \ \ \ \ \ \ }}^{m-j+1}
\;\;\;\;
\\
\left[
\mbox{\scriptsize $
\begin{array}{ccc}
p_0 & &  \\ p_1 & \ddots & \\
\vdots & \ddots & p_0 \\ p_m & & p_1 \\
& \ddots & \vdots \\ & & p_m
\end{array}$} \right. &
\left. \mbox{\scriptsize $\begin{array}{ccc}
q_0 & &  \\ q_1 & \ddots & \\
\vdots & \ddots & q_0 \\ q_n & & q_1 \\
& \ddots & \vdots \\ & & q_n
\end{array}$} \right] \end{array}
\end{equation}
Notice that, for convenience of discussion, we extend the use of
the Sylvester matrices to allow degrees ~$\dg(p) < m$
~and/or ~$\dg(q) < n$ ~in (\ref{Sj}).
For the special case of ~$j=1$ ~and ~$(p,q) \in \cPmn$,
~the matrix ~$S_1(p,q)$ ~is the standard Sylvester matrix in the literature
whose determinant being zero
is frequently used as an equivalent statement for
the existence nontrivial GCD.
~Moreover, the GCD degree can be identified from the nullity of
Sylvester matrices in the following lemma.

\begin{lem} \label{sylthm0}
Let $(p,q)$ be a polynomial pair in $\cPmn$ and $S_j(p,q)$
be the $j$-th Sylvester matrix of $(p,q)$ ~in ~$\mPmn$.
Then the degree of $\GCD(p,q)$ equals to $k$ ~if and only if
\label{'nullity'}
\begin{equation} \label{sylthm}
\nullity{S_j(p,q)} ~=~ \dm\big(\mP_{k-j} \big) ~=~ k-j+1
\end{equation}
for ~$j=1,\ldots,k$ ~and ~$\nullity{S_j(p,q)}\,=\,0$ ~for
~$j = k+1,\ldots,\min\{m,n\}$.
~In particular,
\begin{itemize} \parskip0mm
\item[(i)] $\nullity{S_1(p,q)} = \dg\big(\GCD(p,q)\big)$; and
\item[(ii)] $\nullity{S_k(p,q)} = 1$ with the kernel of
$S_k(p,q)$ ~being spanned by the vector $[\bdw;-\bdv]$ formed by
the cofactors $v$ and $w$ of $(p,q)$.
\end{itemize}
\end{lem}

\prf
Let $v$ and $w$ ~be the GCD cofactors of $(p,q)$.
The kernel of $S_j(p,q)$ ~can be identified from the identity
$(x^i w) \cdot p - (x^i v) \cdot q \;=\; 0$
for $i = 0, 1, \ldots, k-j$ if $j\ge k$. 
\qed

The identity (\ref{sylthm}) in various forms are well known in the 
literature
(see e.g. \cite{emiris-galligo-lombardi,li-zeng-03,rupprecht}),
~while the algorithm in \cite{zeng-dayton,zeng-mr-05} takes advantage of
~$\nullity{S_k(p,q)} \,=\, 1$ ~so that the cofactors ~$v$ ~and ~$w$ can
be solved from the homogeneous linear system
\begin{equation}
S_k(p,q) \left[ \mbox{\scriptsize
$\begin{array}{r} \bdw \\ -\bdv \end{array}$} \right]
~ = ~ \bdo. \label{exvw}
\end{equation}
Then the GCD can be determined via solving the linear system
\begin{equation} \label{sys4u}
C_k(v) \, \bdu ~=~ \bdp  ~~\mbox{and}~~ C_k(w) \, \bdu ~=~ \bdq
\end{equation}
for polynomial ~$u$.

\section{Geometry of GCD and its ill-posedness} \label{sec:geo}

In this section, we study the geometry of the polynomial GCD problem,
the root of its ill-posedness, and the reason why it is not
hypersensitive in a restricted domain in which it becomes numerically
computable.
~The regularization theory that follows later is also derived from
the differentiable manifolds and the stratification structure
formed by the collections of polynomials pairs with common GCD degrees.

Let $(p,q) \in \cPmn$ be a polynomial pair with a particular GCD
$u_* \in \cP_k$ and cofactors $v_*$ and $w_*$.
For any vector $\bdh \in \C^{k+1}$ with
$\bdh^\h \bdu_* = \bt \ne 0$, this GCD triplet $(u_*,v_*,w_*)$ of
$(p,q)$ is the unique solution to the equation
\begin{equation} \label{agcdsys}
\bdf_\bdh(u,v,w) ~=~ [\bt; \,\bdp; \,\bdq] \end{equation}
for $u \in \mP_k$, $v \in \mP_{m-k+1}$, $w \in 
\mP_{n-k+1}$,
where
\begin{equation} \label{fzb}
\bdf_\bdh(u,v,w)  ~~=\;  \mbox{\scriptsize $\left[ \begin{array}{r}
    \bdh^\h \bdu \\ C_k(v) \bdu \\ C_k(w) \bdu
\end{array} \right]$}
\end{equation}
with its Jacobian
\begin{equation}
J_\bdh(u,v,w)
 \;=\;  \mbox{\scriptsize $\left[ \begin{array}{ccc}
\mbox{\normalsize $\bdh^\h$} & & \\
C_k(v) & C_{m-k}(u) & \\
C_k(w) & & C_{n-k} (u) \end{array} \right]$}
\label{jac}
\end{equation}
in which a matrix block such as $C_k(v)$ is the convolution matrix
(\ref{cauchy})
corresponding to the linear transformation $\cL ~:~ g \in \mP_k
\longrightarrow v\cdot g \in \mP_n$.

\begin{lem} \label{jaclem}
Let polynomials $u \in \cP_k$, $v \in \cP_{m-k}$, 
$w \in \cP_{n-k}$ and the vector $\bdh \in \C^{k+1}$
with $\bdh^\h \bdu \ne 0$.
Then the matrix $J_\bdh(u,v,w)$ defined in {\em (\ref{jac})}
is injective if and only if there exists no non-constant polynomial
that divides $u$, $v$ and $w$ simultaneously.
\end{lem}

\prf
Let $a \in \mP_k$, $b \in \mP_{m-k}$ and
$c \in \mP_{n-k}$ be arbitrary polynomials whose
coefficient vectors $\bda \,\in\, \C^{k+1}$,
$\bdb \in \C^{m-k+1}$ and $\bdc \in \C^{n-k+1}$
satisfy
\begin{equation} \label{jnull}
\mbox{\scriptsize $\left[ \begin{array}{ccc}
\mbox{\normalsize $\bdh^\h$} & & \\
C_k(v) & C_{m-k}(u) & \\
C_k(w) & & C_{n-k}(u) \end{array} \right] \;
\left[ \begin{array}{c} -\bda \\ \bdb \\ \bdc
\end{array} \right]$} ~~=~~ \bdo.
\end{equation} \normalsize
The matrix $J_\bdh(u,v,w)$ is injective if $a=b=c=0$.
Equation (\ref{jnull}) is equivalent to
\begin{equation} \label{rp0}
\bdh^\h \bda \;=\; 0, \;\;\; b \cdot u - a\cdot v \;=\; 0, \;\;\;
c \cdot u - a\cdot w \;=\; 0.
\end{equation}
If $u$ and $v$ are co-prime,
then $bu-av=0$ in (\ref{rp0}) implies $a = ru$ for a polynomial $r$.
Then $\dg(a) \le k = \dg(u)$ leads to $r$ being a constant.
Therefore $r=0$ due to $\bdh^\h \bda = r(\bdh^\h \bdu) = 0$.
Namely, $a=0$, which results in $b=c=0$ from (\ref{rp0}).

Assume $u$ and $v$ are not co-prime, namely
$d \in \GCD(u,v)$ is not a constant but $1 \in \GCD(d,w)$.
Write $u = d\cdot u_0$ and $v = d\cdot v_0$.
From $bu-av=0$ in (\ref{rp0}) we have $bu_0-av_0 = 0$ and thus
$u_0$ divides $a$ since $1\in \GCD(u_0,v_0)$, and
$a = s \cdot u_0$ for certain polynomial $s$.
Combining with $cu-aw=0$ yields $c\cdot d - s \cdot w = 0$ and
thus $s = t\cdot d$ for certain polynomial $t$.
Consequently $a = s\cdot u_0 = t\cdot d \cdot u_0 = t\cdot u$.
The polynomial $t$ must be a constant because $\dg(a) \le \dg(u)$.
Moreover $\bdh^\h \bda = t(\bdh^\h \bdu) = 0$ implies $t=0$ and thus
$a = 0$, leading to $b=c=0$ from (\ref{rp0}).
As a result, we have proved $J_\bdh(u,v,w)$ is injective whenever $u$,
$v$ and $w$ have no common non-constant factors.

Assuming there is a non-constant common factor $e$ among
$u$, $v$ and $w$, we now prove $J_\bdh(u,v,w)$ is rank-deficient.
Write $u = e u_0$, $v = e v_0$ and $w = e w_0$.
If $\bdh^\h \bdu_0 = 0$,
then $J_\bdh(u,v,w)[-\bdu_0;\bdv_0;\bdw_0] = \bdo$ by a straightforward
verification, and thus $J_\bdh(u,v,w)$ is rank-deficient.
Next we assume $\bdh^\h \bdu_0 \ne 0$.
Let $a = g u_0$, $b = g v_0$ and $c = g w_0$ for
$g = e - \gamma$ where $\gamma = (\bdh^\h \bdu)/(\bdh^\h \bdu_0)$.
Then $a \in \mP_k$, $b \in \mP_{m-k}$, and
$c \in \mP_{n-k}$ with $bu-av = cu-aw=0$.
The polynomial $g \ne 0$ since $e$ is non-constant, and $a = g u_0 =
d u_0 - \gamma u_0 = u - \gamma u_0$, leading to
$\bdh^\h \bda = \bdh^\h (\bdu - \gamma \bdu_0) = 0$
Therefore $a$, $b$ and $c$ satisfy (\ref{rp0}), implying
$J_\bdh(u,v,w)$ is rank-deficient.
\qed

Lemma~\ref{jaclem} directly leads to the following injectiveness corollary
for the Jacobian (\ref{jac}) at a GCD and cofactors.

\begin{cor} \label{c:jac}
Let $(p,q) \in \cPkmn$.
For almost all $\bdh \in \C^{k+1}$ and $\bt \in \C$, 
there exists a unique GCD
$u_* \in \GCD(p,q)$ with cofactor pair $(v_*, w_*)$ satisfying
{\em (\ref{fzb})}, and the Jacobian $J_\bdh(u,v,w)$ in {\em (\ref{jac})}
is injective at $(u_*,v_*,w_*)$.
\end{cor}

The following GCD Manifold Theorem provides the essential geometric
properties of the GCD problem.
~We adopt a non-abstract notion of a differentiable manifold from 
differential topology: ~A (complex) differential manifold of dimension 
$d$ is a subset that locally resembles the Euclidean space $\C^d$. 
~More specifically, a subset $\Pi \subset \mPmn$ is called a 
differentiable manifold of dimension $d$ if, for every point 
$(p,q) \in \Pi$, there is an open neighborhood $\Delta$ of
$(p,q)$ in $\mPmn$ and a continuously differentiable mapping
$\bdg$ from $\Delta \cap \Pi$ to an open subset $\Lambda$ of $\C^d$ 
with a continuously differentiable inverse 
$\bdg^{-1} : \Lambda \longrightarrow \Delta \cap \Pi$.
The differentiable mapping $\bdg$ is called a {\em local diffeomorphism}
for the manifold $\Pi$, and the {\em codimension} of $\Pi$ is
\[ \codim(\Pi) ~~=~~ \dm(\mPmn) - \dm(\C^d) ~~=~~ m+n+2 - d. \]

\begin{thm}[GCD Manifold Theorem] \label{t:mfld}
With respect to the metric topology induced from the norm $\|\cdot\|$ 
in $\mPmn$, the subset $\cPkmn$ of $\mPmn$ is a differentiable 
manifold of codimension $k$.
Moreover, GCD manifolds $\cPjmn \,\subset\, \overline{\cPkmn}$ if and only
if $j \ge k$, and $\cP^0_{m,n}$ is open dense in $\mPmn$.
\end{thm}

\prf
Let $(p,q)\in \cPkmn$
Then there exist a GCD $\hat{u}$ and its cofactor pair $(\hat{v}$, $\hat{w})$ 
that form the unique solution of the equation (\ref{agcdsys}) for certain 
$\bdh \in \C^{k+1}$ and $\bt \in \C$.
Let the column dimension (and rank) of $J_\bdh(u,v,w)$ be denoted as
\[ l = (k+1)+(m-k+1)+(n-k+1) = m+n-k+3. \]
By Corollary~\ref{c:jac}, there exist $l$ rows of $J_\bdh(u,v,w)$
that are linearly independent.
By the Inverse Function Theorem \cite[Theorem 3.7.3, p.52]{Tay00},
the vector $[\bdu;\bdv;\bdw]$ is locally a continuously differentiable
(vector) function of $l$ entries of the vector $[\bt;\,\bdp;\,\bdq]$.
These $l$ rows must include the first row since otherwise there would be
a contradiction that $(p,q)$ has a unique GCD.
Consequently, the vector $[\bdp;\bdq]$ is a continuously 
differentiable function $\bdg$ of its $l-1$ components in a proper
open domain.
This mapping $\bdg$ is a local diffeomorphism and
the codimension of $\cPkmn$ is thus $(m+n+2)-(l-1) = k$.

The manifold $\cP^0_{m,n}$ is of codimension zero and
thus open in $\mPmn$.
Consequently, the manifold $\cPkmn$ is of positive codimension for
$k > 0 $ and $\mPmn \setminus \cPkmn$ is open dense.
The manifold $\cP^0_{m,n} = \bigcap_{k > 0 }
\big(\mPmn \setminus \cPkmn\big)$ is thus open dense as an intersection of
finitely many open dense subsets.

Let $(p,q) \in \cPkmn$, $\cP^{k'}_{m',n'} \subset
\mPmn$ and
$\inf_{(r,s)\in \cP^{k'}_{m',n'}} \big\|(p,q)-(r,s)\big\| = 0$. 
%
Then there is a sequence
$(p_i,q_i) = (u_i\cdot v_i,\,u_i\cdot w_k) \in \cP^{k'}_{m',n'}$,
$i=1,2,\ldots$ converges to $(p,q)$, where $u_i \in \cP_{k'}$ is
a GCD of $(p_i,q_i)$.
It is clear that $m' = m$ and $n' = n$ since
$(p_i,q_i)$ can not have lower degrees when $i$ is sufficiently large.
Because $\big\{(p_i,q_i)\big\}_{i=1}^\infty$ is bounded, the sequences
$\big\{u_i\big\}_{i=1}^\infty$, $\big\{v_i\big\}_{i=1}^\infty$ and
$\big\{w_i\big\}_{i=1}^\infty$ can be chosen to be bounded and thus
can be assumed as convergent sequences to polynomials
$u \in \mP_{k'}$, $v \in \mP_{m-k'}$
and $w \in \mP_{n-k'}$, respectively.
Consequently, we have $(p,q) = (u\cdot v, \,u\cdot w)$ and thus
$u \in \cP_{k'}$ since otherwise one would have a contradiction in
$\dg(p) < m$ and $\dg(q) < n$.
Therefore $u \in \cP_{k'}$ divides $\GCD(p,q)$.
\qed

The ill-posedness of exact GCD can now be clearly explained: When a
polynomial pair $(f,g) \in \cPkmn$ for $k > 0$ is perturbed,
generically the resulting polynomial pair $(\tilde{f},\tilde{g})$
belongs to $\cPmn^0$ since $\cPkmn$ is dimension
deficient and $\cPmn^0$ is open dense in $\mPmn$.
Consequently the GCD degree drops from $k$ to $0$ discontinuously,
degrading the exact GCD to a constant.
On the other hand, $\cPkmn$ is a differentiable manifold and the
diffeomorphism (\ref{fzb}) has a smooth inverse,
indicating that the GCD is {\em not} discontinuous if the perturbation
is structure-preserving so that $(\tilde{f},\tilde{g})$ remains in
$\cPkmn$.

Theorem \ref{t:mfld} also provides an important geometric property: A small
neighborhood of a polynomial pair $(p,q)$ intersect all the GCD
manifolds $\cPmn^j$ for $j \le k = \dg\big(\GCD(p,q)\big)$.
Furthermore, the residing manifold $\cPkmn$ has a distinct identity
to be given in Lemma~\ref{wolem}.

\section{The notion of numerical GCD} \label{sec:agcd}

We study the numerical GCD for two simultaneous objectives: ~To eliminate
the ill-posedness of the exact GCD and to solve a specific
problem of approximating the GCD that is lost due to data perturbations and
round-off errors.
~The precise problem statement is as follows.

\begin{prob}[The numerical GCD Problem] \label{prob}
Let $(p,q)$ be a given polynomial pair that constitutes the available
data containing a possible perturbation of small magnitude
from an underlying pair $(\hat{p},\hat{q})$.
Find the numerical GCD of $(p,q)$, namely a polynomial $u$ of
degree identical to $\dg\big(\GCD(\hat{p},\hat{q}))$
with an accuracy
\[ \inf_{\hat{u} \in \GCD(\hat{p},\hat{q})} \big\| u-\hat{u}\big\|
~~=~~ O\big(\|(p,q)-(\hat{p},\hat{q})\|\big). \]
\end{prob}

We have regularized ill-posed problems by formulating a ``numerical
solution'' using a ``three-strikes'' principle \cite{ZengNpa,ZengAIF} that
consists of backward nearness, maximum codimension and minimum distance.
~Namely, the numerical solution of the problem is the exact solution of
a nearby problem (backward nearness) that resides in the manifold of the
highest codimension (maximum codimension) and has the minimum distance to the
given data (minimum distance).

We shall introduce the {\em numerical greatest common divisor}
as a well-posed problem to make numerical computation feasible.
~As common in numerical computation, the first and foremost
requirement for computing numerical GCD is its backward accuracy:
~The numerical GCD of a given polynomial pair $(p,q)$
must be the exact GCD of a ``nearby'' pair $(\tilde{p},\tilde{q})$
with $\big\|(p,q)- (\tilde{p},\tilde{q})\big\| < \eps$
for a specified threshold $\eps > 0$.
However, a major distinction here is that
$(\hat{p},\hat{q})$ can not be required as the ``nearest'' pair to $(p,q)$,
as shown in Example~\ref{exam2} below.

\begin{example} \label{exam2}
~Consider the univariate polynomial pair $(p,q)$:
\begin{equation} \label{pq} \left\{\begin{array}{rcl} 
p(x) &=& (x^2-3x+2)\,(x+1.0)+0.01 \\
q(x) &=& (x^2-3x+2)\,(x+1.2)-0.01
\end{array}\right.
\end{equation}
which is small perturbation of magnitude $\sqrt{0.0002} \approx 0.01414$
from a polynomial pair with $\GCD(p,q)=x^2- 3x + 2$ of degree 2.
the nearest polynomial pair with a nontrivial GCD is $(\hat{p}_1,\hat{q}_1)$
where
\begin{eqnarray*} \hat{p}_1 & \approx & 
( x - \mbox{\scriptsize $2.00002$}) (\mbox{\scriptsize $0.9990$}\,x^2 -
\mbox{\scriptsize $0.00255$}\,x - \mbox{\scriptsize $1.0054$}) \\
\hat{q}_1 & \approx & 
( x - \mbox{\scriptsize $2.00002$})
(\mbox{\scriptsize $1.0011$}\,x^2 +
\mbox{\scriptsize $0.2026$}\,x - \mbox{\scriptsize $1.1946$})
\end{eqnarray*}
with distance $0.00168$.
The GCD of $(\hat{p}_1,\hat{q}_1)$ is of degree 1, not a meaningful
approximation to the GCD of degree 2.

In fact, the nearest polynomial pairs with an GCD degree 2 to be
approximately
\[ \big(~
 ( x^2 - \mbox{\scriptsize $3.0001$}\, x + \mbox{\scriptsize $1.9998$})
(\mbox{\scriptsize $1.0017$}\,x + \mbox{\scriptsize $1.0026$}), ~~~
 ( x^2 - \mbox{\scriptsize $3.0001$}\, x + \mbox{\scriptsize $1.9998$})
(\mbox{\scriptsize $0.9982$}\,x +
\mbox{\scriptsize $1.1973$}) ~\big)\]
with larger distance $0.0111$.
If one searches the nearest polynomial pair without a proper constraint,
the actual GCD degree can be misidentified.
\qed
\end{example} 

It is easy to see from Example~\ref{exam2} that, if $(p,q)$ is a polynomial
pair with a nontrivial GCD,  then the pair
$(p,q)$ is {\em closer} to polynomial pairs with GCD of {\em lower}
degrees.
This phenomenon is first reported in \cite{corless-gianni} where it is
suggested to seek the {\em highest degree} for the numerical GCD.
Other than certain non-generic exceptions as we shall see later, this
degree requirement is consistent with a general geometric constraint
for regularizing ill-posed problems:
The numerical GCD must be an exact GCD of a nearby polynomial pair in the
GCD manifold of the {\em highest codimension}.

We shall call $\cPkmn$ the {\em GCD manifold} \,of degree $k$.
A given polynomial pair $(p,q)$ has a distance to each of the
GCD manifolds defined as \label{'thetak'}
\begin{equation} \label{thetak}
 \theta_k(p,q) \;=\; \inf \left\{ \big\|(p,q)-(r,s)\big\| \;\;\Big|\;
(r,s) \in \cPkmn \right\},
\end{equation}
By Theorem~\ref{t:mfld}, those GCD manifolds ~form a
{\em stratification} assuming ~$m \ge n$:
\begin{equation} \label{strat}
\emptyset ~\,=~\, \overline{\cP^{n+1}_{m,n}} ~\,\subsetne~\,
\overline{\cP^n_{m,n}} ~\,\subsetne~\,  \cdots ~\,\subsetne~\,
\overline{\cP^1_{m,n}} ~\,\subsetne~\, \overline{\cP^0_{m,n}}
~\,\equiv~\, \mPmn,
\end{equation}
where $\overline{S}$ denotes the closure of any set $S$.
Consequently, for every $(p,q) \in \cPmn$,
\[ 0 ~=~ \theta_0(p,q) ~\le~ \theta_1(p,q) ~\le~ \cdots ~\le~
\theta_n(p,q). \]
Particularly, for $(p,q) \in \cPkmn$,
\begin{equation} \label{st<}
 0 ~=~ \theta_0(p,q) ~=~ \cdots ~=~ \theta_k(p,q)
~<~ \theta_{k+1}(p,q) ~\le~ \cdots ~\le~ \theta_n(p,q)
\end{equation}
since $k = \dg(\GCD(p,q))$.
The strict inequality in (\ref{st<}) holds because, by Lemma~\ref{sylthm0},
the singular value $\sgmin(S_{k+1}(p,q))$ is strictly positive, while
$\sgmin(S_{k+1}(r,s))=0$ for all polynomial pair
$(r,s) \in \cP^{k+1}_{m,n}$.

\begin{lem} \label{wolem}
Let the pair $(\hat{p},\hat{q}) \in \cPkmn$ and 
$\mathcal{J} \,=\, \big\{ \,j \,\big|\,
\theta_j(\hat{p},\hat{q}) \,=\,0 \big\}$.
Then\linebreak $\mathcal{J} \,=\, \big\{ 0,\,1, \,\ldots\, k\, \}$, namely
$k \,=\, \max \mathcal{J}$.
Furthermore, there exists a $\theta > 0$ such that, from any
$(p,q) \in \cPmn$ with
$\eta = \big\|(p,q)-(\hat{p},\hat{q})\big\| < \theta$, the GCD degree
$k$ of $(\hat{p},\hat{q})$ is identifiable as
\begin{equation} \label{gcdeg}
k ~~=~~
\max{\big\{ j ~\big|~
\theta_j(p,q) \,<\,\eps ~\big\}
}
\end{equation}
for any $\eps$ in the interval $(\eta, \theta)$.
\end{lem}

\prf
A straightforward verification from the GCD Manifold Theorem.
\qed

When the given polynomial pair $(p,q)$ is a small perturbation from
$(\hat{p},\hat{q}) \in \cPkmn$, it can land in any of the
GCD manifold $\cPjmn$ of lower or equal codimension $j \le k$.
However, the underlying GCD degree $k$ distinguishes itself as the
{\em maximum codimension}
\begin{equation} \label{agcddeg}
k ~~\equiv~~ \codim\big(\cPkmn\big) ~~=~~
\max \big\{\,\codim\big(\cPjmn\big)
\;\big|\; \theta_j(p,q) \,<\, \eps \,\big\}
\end{equation}
of all GCD manifolds $\cPjmn$ with
distance $\theta_j(p,q) < \eps$ if the threshold $\eps$ satisfies
\begin{equation} \label{epsineq} \theta_k(p,q) ~<~ \eps ~<~
\theta_{k+1}(p,q), \end{equation}
or the more stringent inequalities $\big\|(p,q)-(\hat{p},\hat{q})\big\|
< \eps < \frac{1}{2} \theta_{k+1}(\hat{p},\hat{q})$.

Revisiting Example~\ref{exam2}, the polynomial pair $(p,q)$ in 
(\ref{pq}) is perturbed from $\cP^2_{4,4}$, which is in the closure 
of $\cP^1_{m,n}$.
By our calculations,
\[ \theta_1(p,q) ~\approx~ 0.00168, ~~~ \theta_2(p,q) ~\approx~ 0.0111
~~~\mbox{and}~~ \theta_3(p,q) ~\approx~ 0.45. \]
The desired GCD manifold ~$\cP^2_{4,4}$ is the one that possesses the highest
codimension 2 and passes through the $\eps$-neighborhood of $(p,q)$ for any
$\eps \in (0.0111,0.225)$.

In this paper, we assume the given polynomial pair $(p,q)$ is a small
perturbation from the underlying pair
$(\hat{p},\hat{q}) \in \cPkmn$, such that
$ \big\|(p,q) - (\hat{p},\hat{q}) \big\|  ~\ll~
\theta_{k+1}(\hat{p},\hat{q})$, and a threshold $\eps$
can be chosen in between and thus (\ref{epsineq}) holds.
If the numerical GCD of $(p,q)$ is of the degree $k$ satisfying
(\ref{agcddeg}), we can recover the underlying GCD degree.
~Furthermore, the minimum distance from $(p,q)$ to the GCD manifold
$\cPkmn$ can be reached at a pair $(\tilde{p},\tilde{q}) \in
\cPkmn$.
We can naturally designate the exact GCD of $(\tilde{p},\tilde{q})$ as 
the numerical GCD of $(p,q)$.

The essential requirements of following numerical GCD definition are first
discovered by Corless, Gianni, Trager and Watt \cite{corless-gianni} in
1995, and formally proposed by Karmarkar and Lakshman in 1996 \cite{kar-Lak96}.

\begin{defn} \label{agcdef}
Let $(p,q) \in \cPmn$ and a threshold $\eps > 0$.
A {\em numerical greatest common divisor} of $(p,q)$ within
$\eps$ is an exact GCD of $(\tilde{p},\tilde{q}) \in \cPkmn$ where $k$
satisfies {\em (\ref{gcdeg})} and
$ \big\| (\tilde{p},\tilde{q}) - (p,q) \big\| \,=\,
\theta_k(p,q)$.
The $\sim$--equivalence class of all numerical GCD's of $(p,q)$
is denoted by $\GCD_\eps(p,q)$.
Namely $\GCD_\eps(p,q) \,=\,\GCD(\tilde{p},\tilde{q})$.
\end{defn}

The formulation of numerical GCD is consistent with
the ``three-strikes principle''
which have been successfully applied to other ill-posed problems
\cite{zeng-mr-05,zeng-li-jcf}.

\begin{itemize} 
\item[] \hspace{-8mm} {\em Backward nearness:} The numerical GCD of a given
polynomial pair $(p,q)$ is the exact GCD of a nearby polynomial pair
$(\tilde{p},\tilde{q})$ within a specified distance $\eps$.
\item[] \hspace{-8mm} {\em Maximum codimension} of the solution manifold:
The nearby pair $(\tilde{p},\tilde{q})$ resides in the highest 
codimension manifold $\cPkmn$ among all the GCD manifolds intersecting
the ``nearness'' $\eps$-neighborhood of the given pair $(p,q)$.
\item[] \hspace{-8mm} {\em Minimum distance} to the solution manifold:
The pair $(\tilde{p},\tilde{q})$ is the nearest point on the manifold
$\cPkmn$ to the given pair $(p,q)$.
\end{itemize} \label{'agcd'}

The numerical GCD defined in Definition~\ref{agcdef} extends the notion
of GCD in the sense that the exact GCD becomes a special case of the 
numerical GCD.
When a pair $(p,q)$ possesses a nontrivial GCD,
the numerical GCD $\GCD_{\eps}(p,q)$ and the exact GCD $\GCD(p,q)$
are identical for all $\eps$ satisfying
$0 < \eps < \theta_{k+1}(p,q)$.

\begin{thm}[GCD Extension Theorem] \label{exthm}
There exists a constant $\theta>0$ associated with  
every polynomial pair $(\hat{p},\hat{q}) \in \cPmn$ possessing an exact 
GCD of degree $k$ such that,
for every $(p,q) \in \cPmn$ that is sufficiently close to 
$(\hat{p},\hat{q})$, there exists a numerical GCD of $(p,q)$ 
within every $\eps \in (0, \theta)$.
This numerical GCD is unique and is of the same degree $k$.
Moreover, 
\[
\lim_{(p,q) \rightarrow (\hat{p},\hat{q})} ~\GCD_\eps(p,q) ~~=~~
\GCD(\hat{p},\hat{q}).
\]
When $(p,q) = (\hat{p},\hat{q})$ in particular, the 
numerical GCD of $(p,q)$ within $\eps \in (0,\theta)$ is identical to 
the exact GCD of $(\hat{p},\hat{q})$.
\end{thm}

\prf
There is a minimum distance $\tau$ from all GCD manifolds having a
positive distance to $(\hat{p},\hat{q})$.
Let $\xi$ be the minimum magnitude of nonzero coefficients of
$(\hat{p},\hat{q})$ and let $\theta = \frac{1}{2}\min\{\tau,\,\xi\}$.
For any $(p,q)$ with $\|(p,q)-(\hat{p},\hat{q})\| < \eps$, the distance
$\theta_k(p,q) < \eps$ and (\ref{gcdeg}) holds, implying
$\dg\big(\GCD_\eps(p,q)\big) = \dg\big(\GCD(\hat{p},\hat{q})\big)$.
The set
\[ \mathcal{S} ~~=~~ \big\{(f,g) \in \mPmn
~\big|~ \|(f,g)-(p,q)\| \le \|(p,q)-(\hat{p},\hat{q})|\big\}
\,\mbox{$\bigcap$}\, \cPkmn
\]
is bounded.
Therefore there exists a convergent sequence $(p_i,q_i) \in \mathcal{S}$
converging to $(p_*,q_*) \in \overline{\cPkmn}$ such
that $\lim_{i\rightarrow \infty} \|(p_i,q_i)-(p,q)\| \,=\,
\theta_k(p,q)$.
Since
\[ \|(p_*,q_*)-(\hat{p},\hat{q})\| ~\le~
\|(p_*,q_*)-(p,q)\| + \|(p,q)-(\hat{p},\hat{q})\| ~\le~ 2\eps ~\le~ \xi,
\]
hence $(p_*,q_*) \in \cPmn$.
If $(p_*,q_*) \not\in \cPkmn$, then $\theta_k(p_*,q_*) = 0$
which lead to\linebreak $\dg\big(\GCD(p_*,q_*)) > k$ by the GCD Manifold Theorem
(Theorem~\ref{t:mfld}), contradicting the choice of $\theta \ge 2\tau$.
Consequently, the distance $\theta_k(p,q)$ is attainable
as $\|(p,q)-(p_*,q_*)\|$, and a $\GCD_\eps(p,q)$ exists.

By Definition~\ref{agcdef}, the equivalence class
$\GCD_\eps(p,q) = \GCD(\tilde{p},\tilde{q})$ where 
$(\tilde{p},\tilde{q}) \in \cPkmn$, and
\[ \big\|(\tilde{p},\tilde{q})-(\hat{p},\hat{q})\big\| \le
\big\|(\tilde{p},\tilde{q})-(p,q)\big\| +
\big\|(p,q)-(\hat{p},\hat{q})\big\|
\le
2\,\big\|(p,q)-(\hat{p},\hat{q})\big\|.
\]
Since $\cPkmn$ is a differentiable manifold and there is a local
diffeomorphism that maps $(u,v,w)$ to $(p,q) \in \cPkmn$ with
$u \in \GCD(\tilde{p},\tilde{q})$, we have
\[ \lim_{(p,q) \rightarrow (\hat{p},\hat{q})} ~\GCD_\eps(p,q) \,=\,
\lim_{(\tilde{p},\tilde{q}) \rightarrow (\hat{p},\hat{q})}
~\GCD(\tilde{p},\tilde{q}) \,=\, \GCD(\hat{p},\hat{q})
\]
and the theorem follows.
\qed

\section{Strong Hadamard well-posedness of numerical GCD} \label{sec:wp}

As introduced by Hadamard, a problem is well-posed (or regular) if its solution
satisfies existence, uniqueness, and certain continuity with respect to data.
~For solving a computational problem accurately using floating point
arithmetic with fixed hardware precision, the continuity must be Lipschitz
so that the Lipschitz constant serves as the finite sensitivity measure,
or otherwise the problem is still incompatible with numerical computation.
~For instance, polynomial roots are continuous with respect to coefficients
regardless of multiplicities.
~However, multiple roots are not Lipschitz continuous and thus infinitely
sensitive to coefficient perturbations, rendering the root-finding problem
extremely difficult until proper regularization is applied
\cite{zeng_multroot,zeng-mr-05,ZengAIF}.
~Consequently, the well-posed problem is often defined in recent literature
as having a finite condition number \cite{DemmelBook}.
~To emphasize the requirement of finite sensitivity, we call the problem
as {\em strongly} well-posed if the continuity is Lipschitz.

We shall establish the strong Hadamard well-posedness of numerical GCD as
formulated in Definition~\ref{agcdef}.
~Particularly, we shall prove a strong well-posedness in Lipschitz continuity.
~To this end, we need the following lemma to prove the regularity of the
numerical GCD.

\begin{lem}
\label{gnlem}
For an open subset $\Omega$ in $\C^n$,
let $\bdf\,:\,\Omega \longrightarrow \C^m$ be analytic
with an injective Jacobian $J(\bdz)$ at every $\bdz\in \Omega$.
Assume ~$\bdf(\Omega)$ ~is a differentiable manifold in ~$\C^n$. 
Then there is an open neighborhood $\Dl$ of every ~$\bdz$ in $\Omega$ and
an open subet ~$\Sigma$ of ~$\bdf(\bdz)$ in ~$\C^m$ such that, for every
$\bdb \in \Sigma$, there is a unique solution $\bdz_* \in \Dl$ to the 
least squares problem
\[ \| \bdf(\bdz_*)-\bdb \|^2 \,=\, 
\min_{\bdy \in \Omega} \|\bdf(\bdy)-\bdb\|^2.
\]
There are also $\sg,\,\gamma > 0$ such that the Gauss-Newton iteration
\begin{equation} \label{gnit}
 \bdz_{k+1} ~=~ \bdz_k - J(\bdz_k)^+\,[\bdf(\bdz_k)-\bdb], ~~~ k = 0, 1, \cdots
\end{equation}
converges to $\bdz_*$ from every initial iterate $\bdz_0 \in \Dl$ 
with
\begin{equation} \label{10250}
 \big\| \bdz_{k+1} - \bdz_* \big\| ~~\leq~~ \sg \,
\big\| \bdz_k - \bdz_* \big\| + \gamma \, \big\| \bdz_k
- \bdz_* \big\|^2
 ~~\leq~~ \mu \, \big\| \bdz_{k} - \bdz_* \big\|
\end{equation}
for $k=0,1,2,\ldots$
with $\mu = \sg + \gamma\, \|\bdz_0-\bdz_*\| < 1$.
Moreover, if $\|\bdf(\bdz_*) - \bdb\|=0$, the convergence rate is 
quadratic with $\sg = 0$.
\end{lem}

\prf
~This is basically a combination of Lemma~2 and Lemma~3 in \cite{ZengAIF}
with a minor variation form the statements of \cite[Lemma 3]{ZengAIF} and
the proof accordingly.
\qed

We now state and prove the following regularity theorem of the numerical
GCD.

\begin{thm}[Numerical GCD Regularity Theorem] \label{regthm}
The numerical GCD\linebreak problem is strongly well-posed.
More specifically,
for every polynomial pair $(\hat{p},\hat{q})$ in $\cPmn$, there is a
neighborhood $\cD$ of $(\hat{p},\hat{q})$ in $\mPmn$ and a
constant $\theta > 0$ such that, for every $(p,q) \in \cD$ and
$\eps$ in the interval $(\|(p,q)-(\hat{p},\hat{q})\|, \,\theta)$,
the following assertions hold:
\begin{itemize} \parskip0mm
\item[(i)] {\em (Existence)} The numerical GCD $\GCD_\eps(p,q)$ 
exists.
\item[(ii)] {\em (Uniqueness)} $\GCD_\eps(p,q)$ is unique in
$\mPmn/\sim$.
\item[(iii)] {\em (Lipschitz continuity)} There is a constant $\alpha >0$
such that, for all $(p_1,q_1)$, $(p_2,q_2) \,\in \cD$, we have
\[  \|(u_1,v_1,w_1) - (u_2,v_2,w_2) \|
\,<\, \alpha\, \big\| (p_1,q_1)-(p_2,q_2) \big\|.
\]
for certain $u_1 \in \GCD_\eps(p_1,q_1)$ and $u_2 \in \GCD_\eps(p_2,q_2)$
with cofactor pairs $(v_1,w_1)$ and $(v_2,w_2)$ respectively.
\end{itemize}
\end{thm}

\prf
The existence of $\GCD_\eps(p,q) = \GCD(\tilde{p},\tilde{q})$ for
$(p,q)$ near $(\hat{p},\hat{q})$ with $(\tilde{p},\tilde{q}) \in
\cPkmn$ is part of the GCD Extension Theorem.
To prove the uniqueness and the Lipschitz continuity,
let $\bdf_\bdh(\cdot,\cdot,\cdot)$ be as in
(\ref{fzb}) along with the Jacobian $J_\bdh(\cdot,\cdot,\cdot)$ as in
(\ref{jac}) with a proper choice of the scaling vector $\bdh$.
Then there is a unique $\hat{u} \in \GCD(\hat{p},\hat{q})$ along with
cofactors $\hat{v}$ and $\hat{w}$ such that
$\bdf_\bdh(\hat{u},\hat{v},\hat{w}) =
[\hat{\beta};\,\hat{\bdp};\,\hat{\bdq}]$
for every scalar $\hat{\bt} > 0$.
Applying Lemma~\ref{gnlem} to $\bdf_\bdh$, there is a neighborhoods
$\Sigma$ of $[\hat{\bt};\,\hat{\bdp};\,\hat{\bdq})$ and
$\Dl$ of
$(\hat{\bdu},\hat{\bdv},\hat{\bdw})$ respectively such that for every
$[\bt;\,\bdp;\,\bdq] \in \Sigma$, there is a unique
$(\tilde{u},\tilde{v},\tilde{w}) \in \Dl$ 
that solves the least squares problem
\[
\|\bdf_\bdh(\tilde{u},\tilde{v},\tilde{w})-[\bt;\bdp;\bdq]\| ~=~
\min_{(u,v,w) \in \mP_k \times \mP_{m-k}\times \in
\mP_{n-k}}
\|\bdf_\bdh(u,v,w)-[\bt;\bdp;\bdq]\|.
\]
Let $(\tilde{p},\tilde{q}) = (\tilde{u}\tilde{v},\,\tilde{u}\tilde{w})$.
Then $(\tilde{p},\tilde{q}) \in \cPmn$ since we can
assume that $\mathcal{D}$ is small so that $\mathcal{D} \subset \cPmn$.
Thus $(u,v,w) \in
\cP_k \times \cP_{m-k}\times \cP_{n-k}$ and
$\dg\big(\GCD(\tilde{p},\tilde{q})\big) \ge \dg(u) = k$.
Since $\cPkmn$ is the GCD manifold of the highest GCD degree near
$(p,q)$ within $\eps$, we have
$\dg\big(\GCD(\tilde{p},\tilde{q})\big) = k$.
Consequently, the uniqueness assertion holds.

Let $[\bt;\,\check{\bdp};\,\check{\bdq}] \,\in\,\Sigma$ and
let $(\check{u},\check{v},\check{w})$ be the least squares solution
to $\bdf_\bdh(\cdot,\cdot,\cdot) = [\bt;\,\check{\bdp};\,\check{\bdq}]$.
Apply one step of the Gauss-Newton iteration on
$\bdf_\bdh(u,v,w) = [\bt;\,\check{\bdp};\,\check{\bdq}]$ from
$(\tilde{u},\tilde{v},\tilde{w})$ and denote
\begin{equation} \label{bdu1}
  [\bdu_1;\,\bdv_1;\,\bdw_1] =
[\tilde{\bdu};\,\tilde{\bdv};\,\tilde{\bdw}] -
J_\bdh(\tilde{u},\tilde{v},\tilde{w})^+
\big( \bdf_\bdh (\tilde{u},\tilde{v},\tilde{w}) -
[\bt;\,\check{\bdp};\,\check{\bdq}] \big).
\end{equation}
Combining (\ref{bdu1}) with
$[\tilde{\bdu};\,\tilde{\bdv};\,\tilde{\bdw}]
\,=\,
[\tilde{\bdu};\,\tilde{\bdv};\,\tilde{\bdw}] -
J_\bdh(\tilde{u},\tilde{v},\tilde{w})^+
\big( \bdf_\bdh (\tilde{u},\tilde{v},\tilde{w}) - [\bt;\,\bdp;\,\bdq] \big)$
yields $\big\|[\bdu_1;\,\bdv_1;\,\bdw_1] -
[\tilde{\bdu};\,\tilde{\bdv};\,\tilde{\bdw}] \big\|
\le
\big\|J_\bdh(\tilde{u},\tilde{v},\tilde{w})^+\big\|
\big\| [\check{\bdp};\,\check{\bdq}] - [\bdp;\,\bdq] \big\|$.
By (\ref{10250}),
\begin{eqnarray*}
\lefteqn{\big\| [\check{\bdu};\,\check{\bdv};\,\check{\bdw}] -
[\tilde{\bdu};\,\tilde{\bdv};\,\tilde{\bdw}] \big\|}\\
& \le &
\big\| [\check{\bdu};\,\check{\bdv};\,\check{\bdw}] -
[\bdu_1;\,\bdv_1;\,\bdw_1] \big\| +
\big\| [\bdu_1;\,\bdv_1;\,\bdw_1]  -
[\tilde{\bdu};\,\tilde{\bdv};\,\tilde{\bdw}] \big\| \\
&  \le &
\mu \big\| [\check{\bdu};\,\check{\bdv};\,\check{\bdw}] -
[\tilde{\bdu};\,\tilde{\bdv};\,\tilde{\bdw}] \big\|
+ \big\|J_\bdh(\tilde{u},\tilde{v},\tilde{w})^+\big\|
\big\|[\bdp;\,\bdq] -[\check{\bdp};\,\check{\bdq}] \big\|
\end{eqnarray*}
Namely
\[
\big\| [\check{\bdu};\,\check{\bdv};\,\check{\bdw}] -
[\tilde{\bdu};\,\tilde{\bdv};\,\tilde{\bdw}] \big\| \le
\frac{ \big\|J_\bdh(\tilde{u},\tilde{v},\tilde{w})^+\big\|}{1-\mu}\,
\big\|[\bdp;\,\bdq] -[\check{\bdp};\,\check{\bdq}] \big\|
\]
where $0 < 1-\mu < 1$
for a sufficiently small $\Sigma$, leading to the Lipschitz continuity.
\qed

Finding numerical GCD not only is a well-posed problem by the
Numerical GCD Regularity Theorem but also solves the problem of 
computing the GCD accurately from perturbed data, as specified in 
Problem~\ref{prob}, by the following Numerical GCD Approximation Theorem.

\begin{cor}[Numerical GCD Approximation Theorem]
~The numerical\linebreak GCD formulated in Definition~\ref{agcdef} solves 
Problem~\ref{prob}.
~More specifically, under the assumptions of Theorem~\ref{regthm},
the numerical GCD ~$\GCD_\eps(p,q)$ ~satisfies the following addition
properties.
\begin{itemize} \parskip0mm
\item[(iv)] {\em (Identical degrees)}
~$\dg\big(\GCD_\eps(p,q)\big) \,=\,
\dg\big(\GCD(\hat{p},\hat{q})\big) \,=\,k$.
\item[(v)] {\em (Convergence)} 
~$\displaystyle \lim_{(p,q) \rightarrow (\hat{p},\hat{q})}
\GCD_\eps(p,q) \,=\, \GCD(\hat{p},\hat{q})$.
\item[(vi)] {\em (Bounded sensitivity)} 
\begin{eqnarray} \label{cond}
\lefteqn{\limsup_{(p,q) \rightarrow (\hat{p},\hat{q})}
\frac{\displaystyle 
\inf_{\mbox{\scriptsize $ u \in \GCD_\eps(p,q),\,
\hat{u} \in \GCD(\hat{p},\hat{q}) 
$}}
\big\|(u,v,w)-(\hat{u},\hat{v},\hat{w})\big\|}{
\|(p,q)-(\hat{p},\hat{q})\|}}  \\
& &  \nonumber \\ & \leq & 
\inf_{\mbox{\scriptsize
$\begin{array}{c} \bdh \in \C^{k+1},\,
\hat{u} \in \GCD(\hat{p},\hat{q}) \\
(\hat{u}\hat{v},\hat{u}\hat{w}) = (\hat{p},\hat{q})
\end{array}$}}
\big\| J_\bdh(\hat{u},\hat{v},\hat{w})^+\big\| 
~~<~~  \infty \nonumber
\end{eqnarray}
where $(v,w)$ and $(\hat{v},\hat{w})$ are cofactor pairs of
$(p,q)$ and ~$(\hat{p},\hat{q})$ respectively.
\end{itemize}
\end{cor}

The above theorem for numerical GCD substantially improves the
similar result in \cite[Proposition~2]{zeng-dayton} and justifies the
definition
\begin{equation} \label{agcdcond}
\kappa_\eps (p,q) ~~=~~
\inf_{\mbox{\scriptsize $\begin{array}{c}
\bdh \in \C^{k+1}, ~~u \in \GCD_\eps(p,q) \\
\|(uv,uw)-(p,q)\|=\theta_k(p,q) \end{array}$}
} \big\| J_\bdh(u,v,w)^+\big\|
\end{equation}
of the {\em numerical GCD condition number} \cite[Definition~2]{zeng-dayton}
of ~$(p,q)$ ~within ~$\eps$.
~We believe the sensitivity measure (\ref{cond}) is optimal.

The condition number $\kappa_\eps(p,q)$ can be estimated as
a by-product of numerical GCD computation.
Upon exit of the Gauss-Newton iteration (\ref{gnit}), the last Jacobian
$J_\bdh(u_i,v_i,w_i)$ is available along with its QR decomposition.
Applying one step of the null vector finder in \cite[p.130]{ZengNpa}
will yield an approximation of the smallest singular value $\sg_{min}$ of
$J_\bdh(u_i,v_i,w_i)$, while $\|J_\bdh(u_i,v_i,w_i)^+\| = 1/\sg_{min}$
can substitute for $\kappa_\eps(p,q)$ as a good estimate.

By Definition~\ref{agcdef}, the GCD Extension Theorem, the Numerical GCD
Regularity Theorem and the Numerical GCD Approximation Theorem, we have 
now established the strong Hadamard well-posedness, and validated the 
so-defined numerical GCD for its intended objective of
solving the numerical GCD Problem as stated in Problem~\ref{prob}.

Lemma \ref{jaclem} provides an insight into the sensitivity of the numerical
GCD by specifying the necessary and sufficient condition for
$J_\bdh(u,v,w)$ to be rank-deficient.
Computing the numerical GCD of $(p,q)$ within $\eps$ is
ill-conditioned if and
only if $J_\bdh(u,v,w)$ is ``nearly'' rank-deficient,
namely $u$, $v$ and $w$ can
be ``nearly'' divisible by a nonconstant polynomial.
Consequently, computing the numerical GCD of $(p,q)$ is not
ill-conditioned even if it is also near a other GCD manifold
as long as the numerical GCD triplet $u$ and cofactors $v$, $w$
do not share an approximate common divisor.

A typical ill-conditioned example can be constructed in the following
example.

\begin{example} ~Consider the following polynomial pair
\begin{equation} \label{ilpair}
\left\{ \begin{array}{lcl}
 p_\dl (x) & = & (x^2-1) \big[ (x-1+\dl) (x^4+1) \big] \\
 q_\dl (x) & = & (x^2-1) \big[ (x-1-\dl) (x^3+2) \big]
\end{array} \right.
\end{equation}
The GCD triplet consists of
\[ u_\dl(x) = x^2-1, ~~ v_\dl(x) = (x-1+\dl) (x^4+1),
~~w_\dl(x) = (x-1-\dl) (x^3+2). \]
For $\dl =0$, there is a common factor $x-1$ among $u_0$, ~$v_0$ and
$w_0$.
Or, $x-1$ ``nearly'' divides all $u_\dl$ ~$v_\dl$ and $w_\dl$.
Consequently, the pair $(p_\dl,q_\dl)$ is ill-conditioned for
$\dl \ll 1$.
Our experiment with {\sc uvGCD} indicates that the condition number
$\kappa_\eps(p_\dl,q_\dl) \approx \frac{1.14}{\dl}$.
\qed
\end{example}

{\bf Remark on formulations of numerical GCD.}
~~In 1985, Sch\"onhage \cite{schonhage} first proposed the {\em quasi-GCD} for
univariate polynomials that needs to satisfy only the backward nearness.
~Sch\"onhage also assumes the given polynomial pair can be
arbitrarily precise even though it is inexact.
~In 1995, Corless, Gianni, Trager and Watt \cite{corless-gianni}
proposed a "highest degree" requirement of GCD in addition to
Sch\"onhage's notion.
~The same paper also suggests minimizing the distance between the given
polynomial pair to the set of pairs with certain GCD degree.
~In 1996/1998 Karmarkar and Lakshman \cite{kar-Lak96,karmarkar-lakshman}
formally defined ``highest degree approximate common divisor problem''
and explicitly included the requirements of backward
nearness, highest degree, and minimum distance.
~It should be noticed that the understanding of numerical GCD
can be significantly different in other works.
~Notably there is another notion of numerical GCD as the
nearest GCD within a certain given degree \cite{kmyz05,KalYanZhi06}.

\section{The initial numerical GCD approximation}

The GCD degree can be identified by the nullity of the Sylvester matrix
as asserted in Lemma~\ref{sylthm0}.
~Likewise, the GCD manifold of maximum codimension specified in the
definition of numerical GCD can be revealed by the {\em numerical} nullity
of the Sylvester matrices.
~The following lemma provides a necessary condition for such a GCD manifold
to be nearby.

\begin{lem} \label{lem:sgmin}
Let $(p,q)$ be a polynomial pair in $\cPmn$ and $\eps >0$.
If the distance $\theta_k(p,q)$ between $(p,q)$ and a GCD manifold
~$\cPkmn$ is less than $\eps$, then
\begin{equation}
\label{sgmin25}
\sg_{-i}\big(S_{j}(p,q)\big)  ~<~ \eps\cdot \sqrt{\max\{m,n\}-j+1}
\end{equation}
for $i=1,2,\ldots, k-j+1$ and $j \le k$,
where $\sg_{-i}\big(S_j(p,q)\big)$ is the $i$-th smallest singular
value of the $j$-th Sylvester matrix for $(p,q)$ in $\mPmn$.
\end{lem}

\prf
Since $\theta_k(p,q) < \eps$, there exists $(r,s) \in \cPkmn$
such that $\big\|(p,q)-(r,s) \big\| < \eps$.
By Lemma~\ref{sylthm0}, singular values $\sg_i\big(S_j(r,s)\big) = 0$ for
$i = 1, 2, \ldots, k-j+1$.
From the linearity
$S_j(p,q) = S_j(r,s) + S_j(p-r, q-s)$
of the Sylvester matrices (\ref{Sj}) and
\cite[Corollary 8.6.2]{golub-vanloan}
\begin{eqnarray*}  \sg_i\big(S_j(p,q)\big) & \le &
\sg_i\big(S_j(r,s)\big) + \big\|S_j(p-r, q-s)\big\|
~~\le~~ \big\|S_j(p-r, q-s)\big\|_F \\
& = & \sqrt{(n-j+1)\|p-r\|^2+(m-j+1)\|q-s\|^2} \\
& <& \eps \cdot \sqrt{\max\{m,n\}-j+1}.
\end{eqnarray*}
\qed

However, inequality (\ref{sgmin25}) does not guarantee the nearness
$\theta_k(p,q) < \eps$, as shown in an example in
\cite{emiris-galligo-lombardi}.
The actual distance $\theta_k(p,q)$ can nonetheless be calculated
during the subsequent computation to ensure finding the numerical GCD
accurately.

\begin{lem} \label{lem:vw}
For a given $(p,q) \in \cPmn$ and $\eps > 0$, let
$(\tilde{p},\tilde{q}) \in \cPkmn$ be the polynomial pair that defines
$\GCD_{\eps}(p,q) = \GCD(\tilde{p},\tilde{q})$ containing $\tilde{u}$ 
with cofactors $\tilde{v}$ and $\tilde{w}$.
If $[\bdw;-\bdv]
\in \C^{n-k+1} \times \C^{m-k+1}$
is the singular vector of ~$S_k(p,q)$ with 
$\big\|S_k(p,q)[\bdw;-\bdv]\big\| = \sgmin\big(S_k(p,q)\big)$, then
$\sgn{2}\big(S_{k}(\tilde{p},\tilde{q})\big) \ne 0$ and the distance
\begin{equation}
\label{cferr}
\dist{
\spn\left\{\left[
\mbox{$\begin{array}{c} \bdw \\ -\bdv \end{array}$}
\right]\right\},
\spn\left\{\left[
\mbox{$\begin{array}{c} \tilde\bdw \\ -\tilde\bdv \end{array}$}
\right] \right\} }
 <  \frac{ 2\,\eps\,\sqrt{\max\{m,n\}-k+1} }{
\sgn{2}\big(S_{k}(\tilde{p},\tilde{q})\big)}.
\end{equation}
\end{lem}

\prf
From Lemma \ref{sylthm0}, we have
$\sgn{2}\big(S_{k}(\tilde{p},\tilde{q})\big) \ne 0$.
Consider the singular value decomposition
\( S_k(\tilde{p},\tilde{q}) = \big[\tilde{U},\tilde{\bdz}\big]
\mbox{\scriptsize
$\left[ \begin{array}{cc} \tilde\Sigma & \\ & 0 \end{array} \right]$}
\big[ \tilde{V}, \tilde{\bdy} \big]^*
\)
and let $\bdy = [\bdw; -\bdv]$.
~We have
\begin{eqnarray*}
\big\| S_k(\tilde{p},\tilde{q}) \bdy \big\| & = &
\big\|\tilde\Sigma \tilde{V}^\h \bdy\big\| ~\ge~
\sgn{2}\big(S_{k}(\tilde{p},\tilde{q})\big) \big\|\tilde{V}^\h\bdy\big\|
~~~~~~~~\mbox{and} \\
\big\| S_k(\tilde{p},\tilde{q}) \bdy \big\| & \le &
\big\| S_k(p,q) \bdy \big\| + \big\|S_k(\tilde{p},\tilde{q}) - S_k(p,q)\big\|\,
\|\bdy\|\\
&<& 2\eps \sqrt{\max\{m,n\}-k+1}.
\end{eqnarray*}
Therefore, the inequality (\ref{cferr}) follows from the identity
\cite[Theorem 2.6.1]{golub-vanloan}
\[
\dist{
\spn\left\{\left[
\mbox{\scriptsize $\begin{array}{c} \bdw \\ -\bdv \end{array}$}
\right]\right\},
~\spn\left\{\left[
\mbox{\scriptsize $\begin{array}{c} \tilde\bdw \\ -\tilde\bdv \end{array}$}
\right] \right\} } ~=~\big\|\tilde{V}^\h\bdy\big\|.
\]
\qed

Lemma \ref{lem:sgmin} provides mechanisms for identifying the
numerical GCD degree and numerical cofactor pair.
~When inequality (\ref{sgmin25}) holds then it is {\em possible} to have an
numerical GCD degree ~$k$, and ~$(v,w)$ ~can be extracted
from the right singular vector.
~The smallest singular value and the corresponding right singular vector can
be computed accurately and efficiently using a numerical rank-revealing
iteration \cite{li-zeng-03,zeng-mr-05} in the following lemma.

\begin{lem} \label{svcvrt}
~Under the assumptions of Lemma {\em \ref{lem:sgmin}},
assume the inequality
\[ \sgn{2}\big(S_k(p,q)\big) ~~>~~ 2 \eps \sqrt{\max\{m,n\}-k+1}
\]
holds and ~$Q\cdot R$ ~is the QR decomposition 
{\em \cite[\S5.2]{golub-vanloan}} of ~$S_k (p,q)$.
~Then, for almost all initial vector ~$\bdz_0$ ~of proper dimension,
the following iteration
\begin{equation} \label{rrit}
\left\{ \begin{array}{l}
\mbox{Solve}~~ R^\h \bdy_j = \bdz_{j-1} ~~\mbox{by forward substitution} \\
\mbox{Solve}~~ R\, \bdz_j = \bdy_j ~~\mbox{by backward substitution} \\
\mbox{Normalize}~~ \bdz_j,  ~~\mbox{for}~~ j = 1, 2, \ldots
\end{array} \right.
\end{equation}
generates a sequence of unit vectors $\bdz_j$, $j=1,2,\ldots$ converging
to ~$\bdz_*$ ~and
\begin{equation} \label{rrit1}
 \|S_k(p,q)\, \bdz_*\| ~=~ \|R\, \bdz_*\| ~=~
\sgn{1}\big(S_k(p,q)\big)
\end{equation}
at convergence rate
\begin{equation} \label{rritrate}
 \|\bdz_j - \bdz_*\| ~\le~
\big[\sgn{1}\big(S_k(p,q)\big)\big/
\sgn{2}\big(S_k(p,q)\big)\big]^{2j} \|\bdz_0 - \bdz_*\|
\end{equation}
\end{lem}

\prf
~From $\sgn{2}\big(S_k(p,q)\big) > 2 \eps \sqrt{\max\{m,n\}-k+1}$ and
Lemma \ref{lem:sgmin}, we have
\begin{eqnarray*}
\sgn{1}\big(S_k(p,q)\big) & \le & \eps
\sqrt{\max\{m,n\}-k+1} ~~\mbox{and}~~  \\
\sgn{2}\big(S_k(p,q)\big) & \ge &
\sgn{2}\big(S_k(p,q)\big) - \eps \sqrt{\max\{m,n\}-k+1}\\
 &>& \eps \sqrt{\max\{m,n\}-k+1}
\end{eqnarray*}
and thus ~$\sgn{1}\big(S_k(p,q)\big)\big/
\sgn{2}\big(S_k(p,q)\big) ~<~ 1$.
~The assertions of the lemma then follows from \cite[Lemma 2.6]{zeng-mr-05}.
\qed

Equations in (\ref{rrit1}) implies ~$\bdz_*$ ~is the vector
~$[\bdw; -\bdv]$ ~in Lemma~\ref{lem:vw} containing the
coefficients of the numerical cofactors ~$v$ ~and ~$w$.
~The next lemma provides an error estimate for the initial approximation
~$u$ ~of the numerical GCD from solving the least squares solution to
system ~(\ref{sys4u}).

\begin{lem} \label{lem:uerr}
Under the assumptions of Lemma {\em \ref{lem:sgmin}} and {\em \ref{lem:vw}} 
with the same notations along with $\mu\, \eps$ denoting the
right hand side of {\em (\ref{cferr})}, let 
\[ \xi = \mbox{\scriptsize $
\left\| \left[ \begin{array}{c} C_k(\tilde{v}) \\ C_k(\tilde{w}) \end{array} 
\right] \right\|$},  ~~~~~\tau ~=~
\xi \,\mbox{\scriptsize $
\bigg\| \left[ \begin{array}{c} C_k(\tilde{v}) \\ C_k(\tilde{w}) \end{array}
\right]^+\bigg\|$}, 
\]
and $\bdz \,=\, \bdu$~ be the least squares solution to
\begin{equation} \label{sys4uhat}
\mbox{\scriptsize $
\left[ \begin{array}{c} C_k(v) \\ C_k(w) \end{array}
\right]$} \,\bdz \;=\;
 \mbox{\scriptsize $
\left[ \begin{array}{c} \bdp \\ \bdq \end{array} \right]$}.
\end{equation}
If ~$\eta \,=\, \mu \tau \sqrt{k+1} \,\eps \;<\; 1$,
~then there is an $\al \in \C \setminus \{0\}$ such that
\begin{equation} \label{uhaterr}
\Big\|\tilde{u} - \al\, u\, \Big\| \;\equiv\;
\Big\|\tilde\bdu - \al\,\bdu \Big\| \;\le\;
\mbox{$\frac{\tau}{1-\eta} \left[ \sqrt{k+1}\,\|\tilde{u}\|
 \mu +\frac{1}{\xi} \right]$}~\eps
\end{equation}
\end{lem}

\prf
Let $A \,=\,
\mbox{\scriptsize $ \left[ \begin{array}{c} C_k(\tilde{v}) \\ C_k(\tilde{w}) 
\end{array} \right]$}$~ and ~$\bdb \,=\,
\mbox{\scriptsize $ \left[ \begin{array}{c} \tilde\bdp \\ \tilde\bdq 
\end{array} \right]$}$.
The overdetermined linear system $A \bdz \,=\, \bdb$ has a conventional
solution $\bdz = \tilde{\bdu}$.
Due to (\ref{cferr}), there is a $\gamma \in \C \setminus \{0\}$ ~such that 
$\|\gamma (v,w) - (\tilde{v},\tilde{w})\| \le \mu\,\eps 
\| (\tilde{v},\tilde{w})\|$. 
Rewrite the linear system (\ref{sys4uhat}) as
\[ \mbox{\scriptsize $ \left[ \begin{array}{c} C_k(\gamma v) \\ C_k(\gamma w) 
\end{array} \right]$} (\bdz/\gamma) = 
\mbox{\scriptsize $ \left[ \begin{array}{c} \bdp \\ \bdq 
\end{array} \right]$}
\]
that can be considered as the perturbed system
~$(A + \dl A)(\bdz + \dl \bdz) \,=\, \bdb + \dl \bdb$
where
\begin{eqnarray*}
\lefteqn{\|\dl A \| ~ = ~
\left\| \mbox{\scriptsize $\left[ \begin{array}{c}
C_k(\gamma v - \tilde{v}) \\
C_k(\gamma w - \tilde{w}) \end{array}\right]$} \right\| \;\le\;
\left\| \mbox{\scriptsize $\left[ \begin{array}{c} C_k(\gamma v-\tilde{v}) \\
C_k(\gamma w-\tilde{w}) \end{array}\right]$} \right\|_F }\\
& &  ~~~~~~~~\le\;
\sqrt{k+1} \cdot \big\|\gamma (v,w)-(\tilde{v},\tilde{w})\big\| \;=\;
\sqrt{k+1}
\,\mu \,\eps \|(\tilde{v},\tilde{w})\|, \\
&&\|A\|  ~\ge ~ \frac{1}{\sqrt{k+1}} \,\|A\|_F \;=\;
\frac{1}{\sqrt{k+1}} \sqrt{k+1} \,
\|(\tilde{v},\tilde{w})\| \;=\; \|(\tilde{v},\tilde{w})\|,\\ 
&&\|\bdb\| ~ = ~ \|(\tilde{p},\tilde{q})\|, ~~~~~
\|\dl \bdb \|   \;=\;
 \mbox{\scriptsize $\left\| \left[
\begin{array}{c} \bdp - \tilde\bdp \\ \bdq- \tilde\bdq \end{array} \right]
\right\|$} \;\le\; \eps, ~~~~ \tau \;=\; \|A\|\left\|A^+\right\|,
\end{eqnarray*}
Then inequality (\ref{uhaterr}) follows from Theorem 1.4.6 and Remark 1.4.1 
in \cite[pp. 30-31]{bjorck96}, residual ~$\|A\bdu - \bdb\| \,=\,0$ and
$\al = 1/\gamma$.
\qed

Lemma~\ref{lem:sgmin} and Lemma~\ref{lem:uerr} lead to the following
lemma that ensures the initial approximation of the numerical GCD
and cofactors to be sufficiently accurate if the perturbation to the
polynomial pair ~$(p,q)$ ~is small, satisfying the local convergence
condition of the Gauss-Newton iteration given in Lemma~\ref{gnconv}.

\begin{lem} \label{prop:tricont}
Let ~$(\tilde{p},\tilde{q}) \in \cPkmn$ and $\tilde{u} \in 
\GCD(\tilde{p},\tilde{q})$ with cofactor pair $(\tilde{v},\tilde{w})$.
Then for any ~$\dl > 0$, there is an $\eta > 0$ ~such that for all
$(p,q) \in \cPmn$ with distance
$\big\|(\tilde{p},\tilde{q})-(p,q)\|\,<\,\eta$, the inequality
$\big\|(\tilde{u},\tilde{v},\tilde{w})-(\frac{1}{\gamma}u, \gamma v, 
\gamma w)\big\|
\,<\, \dl$ holds for certain $\gamma \in \C \setminus \{0\}$
where $(u,v,w)$ is defined in Lemma
\ref{lem:vw} and Lemma \ref{lem:uerr} corresponding to
$(p,q)$.
\end{lem}

\prf ~A straightforward verification using Lemmas \ref{lem:sgmin}, 
\ref{lem:vw} and \ref{lem:uerr}. \qed

In summary, to calculate the numerical GCD of a given polynomial pair
$(p,q)$ within a prescribed threshold ~$\eps$,
we first identify the numerical GCD
degree $k$.
Lemma \ref{lem:sgmin} suggests that we can calculate the
smallest singular value
$\sgmin\big(S_j(p,q)\big)$ for $j$ ~decreasing from $\min \{m,n\} = n$
until $\sgmin\big(S_j(p,q)\big) \,<\, \eps\,\sqrt{m-k+1}$
and set $k\,=\,j$.
After $k = \dg\big(\GCD_\eps(p,q)\big)$ is determined, the corresponding
singular vector of $S_{k}(p,q)$ provides an approximation
$(v,w)$ to $(\tilde{v},\tilde{w})$ with error bound (\ref{cferr}).
An approximation $u$
to $\tilde{u} \in \GCD_{\eps}(\tilde{p},\tilde{q})$ 
is obtained by
solving the overdetermined linear system (\ref{sys4uhat}) for the least
squares solution ~$\bdz \,=\, \bdu$~ with error bound (\ref{uhaterr}).
The triplet $(u,v,w)$ will be taken as an initial
iterate for the Gauss-Newton iteration (\ref{gnit3}) for verification and
refinement.

\section{Sensitivity of numerical GCD computation via Sylvester matrices}

The sensitivity of the triplet
~$(u,v,w)$ ~in Lemma \ref{lem:vw} and
Lemma \ref{lem:uerr} can be measured by the reciprocal of
~$\sgn{2}\big(S_{k}(\tilde{p},\tilde{q})\big)$,
as indicated by inequalities in  (\ref{cferr}) and (\ref{uhaterr}).
~In other words, computing the triplet ~$(u,v,w)$ ~by
iteration (\ref{rrit}) ~in combination with solving the linear system
(\ref{sys4uhat}) ~is ill-conditioned whenever 
~$\sgn{2}\big(S_{k}(\tilde{p},\tilde{q})\big)$ ~is tiny.
~Such ill-condition is certain to occur when the pair
~$(p,q) \in \cPkmn$ ~is also near
another GCD manifold ~$\cPjmn$ ~of {\em higher} GCD degree
~$j$.
~We can actually estimate the magnitude of 
~$\sgn{2}\big(S_k(\tilde{p},\tilde{q})\big)$ ~as follows.

Let ~$(\hat{p},\hat{q}) \in \cPjmn$ ~with degree ~$j > k$
~and the distance ~$\big\|(\tilde{p},\tilde{q})-(\hat{p},\hat{q})\big\| = 
\dl$ ~being small.
~By Lemma~\ref{sylthm0}, ~$\nullity{S_k(\hat{p},\hat{q})} = j-k+1 > 2$.
~Similar to the proof of Lemma~\ref{lem:sgmin},
\begin{eqnarray*}
 \sgn{2}\big(S_\bdk(\tilde{p},\tilde{q})\big) & \le & 
\sgn{2}\big(S_\bdk(\hat{p},\hat{q})\big) +
\dl \cdot \sqrt{\max\{m,n\}-k+1} \\
& = & 
\dl \cdot \sqrt{\max\{m,n\}-k+1}.
\end{eqnarray*}
Roughly speaking, the error of the numerical GCD triplet
~$(u,v,w)$ ~computed as in Lemma~\ref{lem:vw} and
Lemma~\ref{lem:uerr} is inversely proportional to the
distance between the polynomial pair ~$(\tilde{p},\tilde{q})$ ~and the 
nearest GCD manifold of higher codimension.
~The following is a typical example in which the polynomial pair is sensitive
for computing the initial numerical GCD approximation but well-conditioned
if it is measured by the GCD condition number.

\begin{example}
Consider the polynomial pair $p_\mu = u\cdot v_\mu$ and
$q = u\cdot w$ where
\[ u(x) ~=~ x^2+1, ~~v_\mu(x) = (x - 1 +\mu)(x^4+1),
~~w(x) = (x - 1)(x^3-2)  \]
Clearly, $(p_\mu,q) \in \cP^2_{7,6}$ with
$\GCD(p_\mu, q) = u$ for all $\mu \ne 0$, but
$(p_\mu,q)$ is near $\cP^3_{7,6}$ when $\mu$ is small.
In fact, the distance $\theta_3(p_\mu,q)$ between $(p_\mu,q)$
and $\cP^3_{7,6}$ is bounded by $\|r_\mu\|=2\mu$ for
$r_\mu = \mu(x^2+1)(x^4+1)$.
While $\sgn{2}\big(S_2(p_\mu,q)\big) > 0$ by Lemma \ref{sylthm} but
the nullity of $S_2(p_\mu-r_\mu,q)$ is at least 2 since
\begin{eqnarray*}  (p_\mu-r_\mu)\cdot (x^3-1)-q \cdot (x^4+1) &=& 0, 
~~\mbox{and}~~ \\
 (p_\mu-r_\mu)\cdot (x-1)(x^3-1)-q \cdot (x-1)(x^4+1) &=& 0. 
\end{eqnarray*}
Hence $\sgn{2}\big(S_2(p_\mu -r_\mu,q)\big) = 0$.
As a result, 
\[ \sgn{2}\big(S_2(p_\mu,q)\big) \le \|C_{6-2}(r_\mu)\|_F =
2\sqrt{5} \mu, 
\]
and it is sensitive to compute the numerical GCD
solely relying on Lemma~\ref{lem:vw} and Lemma~\ref{lem:uerr}.
For instance, let $\mu = 10^{-12}$.
A straightforward computation of $(\tilde{u},\tilde{v},\tilde{w})$ in
Matlab by Lemma~\ref{lem:vw} and Lemma~\ref{lem:uerr} results only three
to four digits accuracy:
\begin{eqnarray*} \tilde{u}(x)  & \approx &  x^2 + 0.99992 \\
 \tilde{v}(x) &\approx& x^5 - 0.9995x^4 + x - 0.9995 \\
\tilde{w}(x)  &\approx&   x^4 - 0.9995x^3 - 2x + 1.9991
\end{eqnarray*}
It may seem to be a surprise that computing $\GCD_{\eps}(p_\mu,q)$
is {\em not} ill-conditioned even if $\mu \ll 1$.
The numerical GCD condition number is nearly a constant of moderate magnitude
($\approx 3.55$) for varying $\mu$.
Even for $\mu = 10^{-12}$, our software {\sc uvGCD} still calculates
the numerical GCD with an accuracy around machine precision
($\approx 2.2\times
10^{-16})$.
The reason for such a healthy numerical condition is revealed in
Lemma~\ref{jaclem}: ~Even though
polynomials $v_\mu$ and $w_\mu$ are close to having a nontrivial common
factor $(x-1)$, the GCD triplet members $u,v_\mu$, and $w$ as a whole
are not.
\qed
\end{example} 

The sensitivity analysis the example above show that, to ensure 
accuracy, it is essential to refine the numerical GCD after obtaining an
initial approximation to the numerical GCD and cofactors.
~Such refinement can be carried out by the Gauss-Newton iteration
that is to be discussed in the next section.

\section{Minimizing the distance to a GCD manifold} \label{s:min}

For a given polynomial pair $(p,q) \in \cPmn$ with the
degree $k = \dg\big(\GCD_\eps(p,q)\big)$ of the numerical GCD being
calculated from Lemma~\ref{lem:sgmin},
finding its numerical GCD and cofactors becomes the problem of minimizing
the distance from $(p,q)$ to the GCD manifold $\cPkmn$:
\[ \big\|(p,q)-(u\cdot v,~u\cdot w) \big\| ~~=~~
\min_{(r,s)\in\cPkmn}\Big\|(p,q)-(r,s)\Big\|,
\]
where $\dg(u) = k$, $\dg(v) = m-k$ and
$\dg(w)=n-k$.
Naturally, this minimization leads to the least squares
problem for the quadratic system consists of
\begin{equation} \label{vupwuq}
~C_k(v) \,\bdu ~=~ \bdp,  ~~~C_k(w)\,\bdu ~=~ \bdq
\end{equation}
which are the vector form of $u\cdot v = p$ and $u\cdot w = q$
respectively.
However, the system (\ref{vupwuq}) is not regular since the least squares
solutions are not isolated.
Any solution $(u,v,w)$ can be arbitrarily scaled to
$(\al u, ~v/\al, ~w/\al)$.
A simple auxiliary equation $\bdh^\h\bdu = \bt$ takes away this dimension
of the solution and ensures the Jacobian to be injective.

We minimize the distance from a point $(uv,uw)$ in the GCD manifold to
the give polynomial pair $(p,q)$ by solving the system
$\bdf_\bdh(u,v,w) \,=\, [\bt;\,\bdp;\,\bdq]$
as in (\ref{agcdsys})
for its least squares solution, where the function
$\bdf_\bdh(u,v,w)$ is defined in (\ref{fzb}).
The Gauss-Newton iteration (\ref{gnit}) for finding
$u_* \in \GCD_{\eps}(p,q)$ and cofactors becomes
\begin{equation}
\left[ \begin{array}{c} 
\bdu_{j+1} \\ \bdv_{j+1} \\ \bdw_{j+1} 
\end{array} \right] = 
\left[ \begin{array}{c} 
\bdu_j \\ \bdv_j \\ \bdw_j
\end{array} \right] 
 - J_\bdh(u_j,v_j,w_j)^+
\left( \bdf_\bdh(u_j,v_j,w_j) - 
\left[ \begin{array}{c} \bt \\ \bdp \\ \bdq \end{array} \right] 
\right)\label{gnit3} 
\end{equation}
for $j=0,1,\ldots$ where $J_\bdh(\cdot,\cdot,\cdot)$ is the Jacobian of
$\bdf_\bdh(\cdot,\cdot,\cdot)$ given in (\ref{jac}).
Lemma \ref{jaclem} ensures this iteration to be locally convergent
for finding the least squares solution $(\hat{u},\hat{v},\hat{w})$
to the system $\bdf_\bdh(u,v,w) = [\bt;\,\bdp;\,\bdq]$.

\begin{lem} \label{gnconv}
For $(p,q) \,\in\, \cPmn$ with numerical GCD $u_* \in \GCD_\eps(p,q)$
and cofactor pair $(v_*,w_*)$, let $\bdh \in \C^{k+1}$ and
$\bt = \bdh^\h \bdu_* \ne 0$.
Define $\bdf_\bdh(\cdot,\cdot,\cdot)$ as in
{\em (\ref{fzb})}.
There is a $\rho>0$ such that, if
$\big\| (u_*v_*,u_*w_*) - (p,q) \big\| < \rho$, 
there exists a $\mu>0$ and the Gauss-Newton iteration {\em (\ref{gnit3})}
converges to $(u_*,v_*,w_*)$ from any initial iterate
$(u_0,v_0,w_0)$ satisfying $\bdh^\h \bdu_0 = \bt$
and $\big\| (u_0,v_0,w_0) - (u_*,v_*,w_*) \big\| < \mu$.
\end{lem}

\prf
~The proof is a straightforward verification using Lemma \ref{gnlem}.
\qed.

\section{The two-staged univariate numerical GCD algorithm}

Based on the general analysis in previous sections, we present the algorithm
originally proposed in \cite{zeng-mr-05}
for computing the numerical GCD triplet $(u,v,w)$ of a given polynomial pair
$(p,q) \in \cPmn$ within a given tolerance
$\eps$ of backward error $\big\|(p,q) - (uv,uw) \big\|$.
The algorithm consists of two stages.
At opening stage, we calculate the degree $k$ of the numerical GCD and an
initial approximation $(u_0,v_0,w_0)$ to $(u,v,w)$.
Then the Gauss-Newton iteration is applied to generate a sequence
$(u_j,v_j,w_j)$ such that $(p_j,q_j) = (u_j v_j, u_j w_j) \in
\cPkmn$ converges to $(\tilde{p},\tilde{q})$ that is the
nearest point on the manifold $\cPkmn$ to the given pair
$(p,q)$.

For simplicity, we assume polynomials $p$ and $q$ are arranged such
that $(p,q)\in \cPmn$ with $m \ge n$ in this section.

\subsection{The numerical GCD degree and the initial GCD approximation}
\label{sec:ini}

Let polynomials $p$ and $q$ be given along with backward error
tolerance $\eps$.
From Lemma \ref{lem:sgmin}, there are no numerical GCD's of degree $j$
within $\eps$ when the smallest singular value
\[ \sgmin\big(S_j(p,q)\big) ~>~
\eps \sqrt{\max\{m,n\}-j+1} ~=~ \eps \,\sqrt{m - j + 1}.
\]
The first stage of numerical GCD computation is to calculated
$\sgmin\big(S_j(p,q)\big)$ for $j$ decreasing from $\min\{m,n\}=n$
through $n-1,n-2,\ldots$ to exclude the possibility of numerical GCD of
those degrees.
The process {\em tentatively} stops when
$\sgmin\big(S_{k}(p,q)\big) \leq \eps \sqrt{m-k+1}$,
pending certification at the refinement stage.

The full singular value decompositions of $S_j(p,q)$'s are unnecessary.
Only the smallest singular value and the associated right singular
vector are needed for each $S_j(p,q)$.
The iteration (\ref{rrit}) is specifically designed for our
purpose here.
It requires the QR decomposition
$Q_jR_j = S_j(p,q)$.
The straightforward computation of each $R_j$
requires $O(j^2 n)$ floating point operations (flops), and the whole
process may require $O(n^4)$ which is unnecessarily expansive.
A successive QR updating strategy as follows substantially reduces the
total flops to $O(n^3)$.

We first calculate the QR decomposition of
$S_n(p,q) = \big[\,C_0(p) \big| C_{m-n}(q)\,\big]$:
\[ S_n(p,q) = Q_n R_n, \;\;\;\;
R_n = \mbox{\scriptsize $\left[ \begin{array}{ccc} * & \ldots & *
\\ 0 & \ddots & \vdots \\ 0 & \ddots & * \\ \vdots & & 0 \\
\vdots & &  \vdots \\ 0 & \ldots & 0 \end{array} \right]$
}_{m \times (m-n+2)}
\]
where each ``$*$'' represents an entry that is potentially nonzero.
When $S_j(p,q)$ is formed for $j=n,n-1,\ldots$, the Sylvester matrix
$S_{j-1}(p,q)$ is
constructed by appending a zero row at the bottom of $S_{j}(p,q)$ followed
by inserting two columns $[\bdo;\,\bdp]$ and $[\bdo; \,\bdq]$.
By a proper column permutation $P_{j-1}$, these two columns are
shifted to the right side as the last two columns of
$S_{j-1}(p,q)P_{j-1}$.
Let initial column permutation $P_n = I$.
Then the QR decomposition $S_n(p,q)P_n = Q_n R_n$ is available.
If the QR decomposition $S_{j}(p,q)P_j = Q_{j}R_{j}$ is available for
$j=n,n-1,\ldots$, the expansion from $S_j(p,q)P_j$ to
$S_{j-1}(p,q)P_{j-1}$ can be illustrated as
\[
\mbox{\raisebox{0.5in}{$S_j(p,q)P_j ~~=~~$}}
\mbox{\epsfig{figure=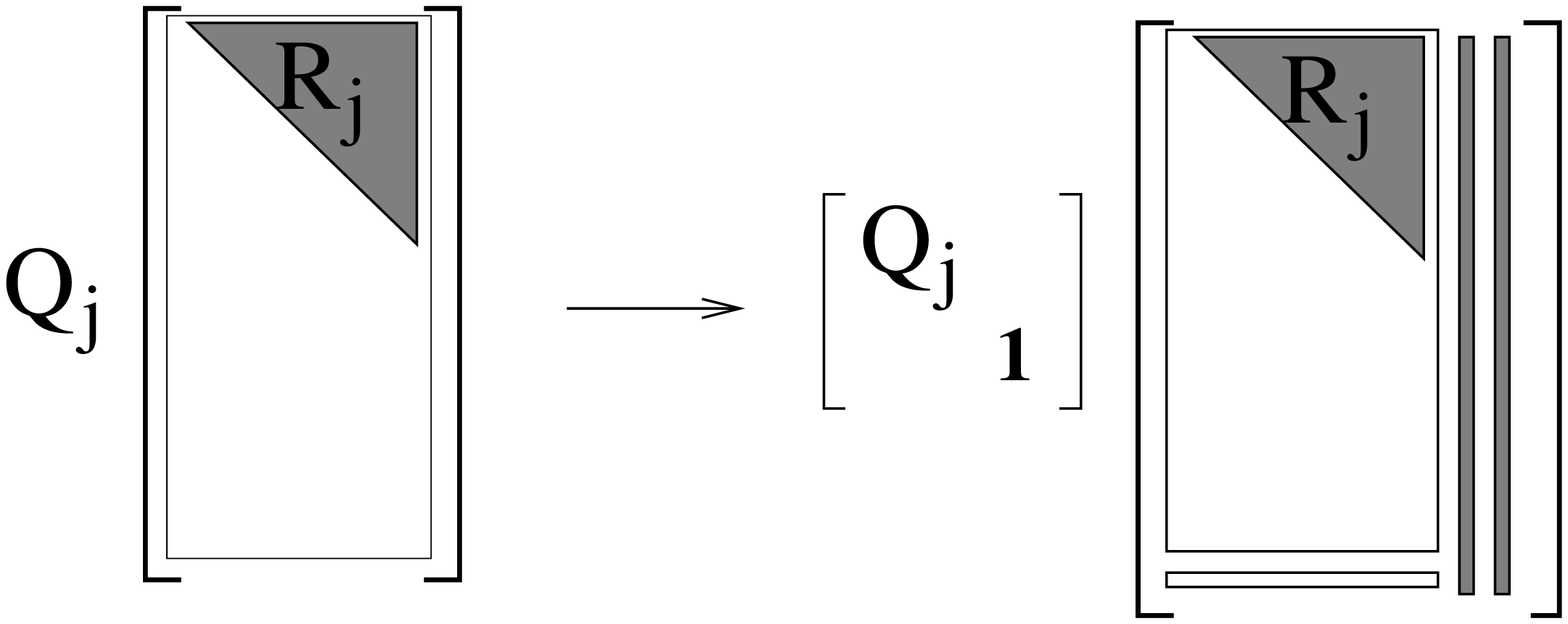,height=1.0in,width=3.0in}}
\mbox{\raisebox{0.5in}{$~~=~~ S_{j-1}(p,q)P_{j-1}$}}
\]
Accordingly, we update $Q_{j}R_{j}$ to $Q_{j-1}R_{j-1}$
by eliminating the lower triangular entries using the Householder
transformation and obtain
\begin{equation}
\mbox{$S_{j-1}(p,q)P_{j-1} \;\;  = $ \ \ \ }
\mbox{\raisebox{-0.5in}{
\mbox{\epsfig{figure=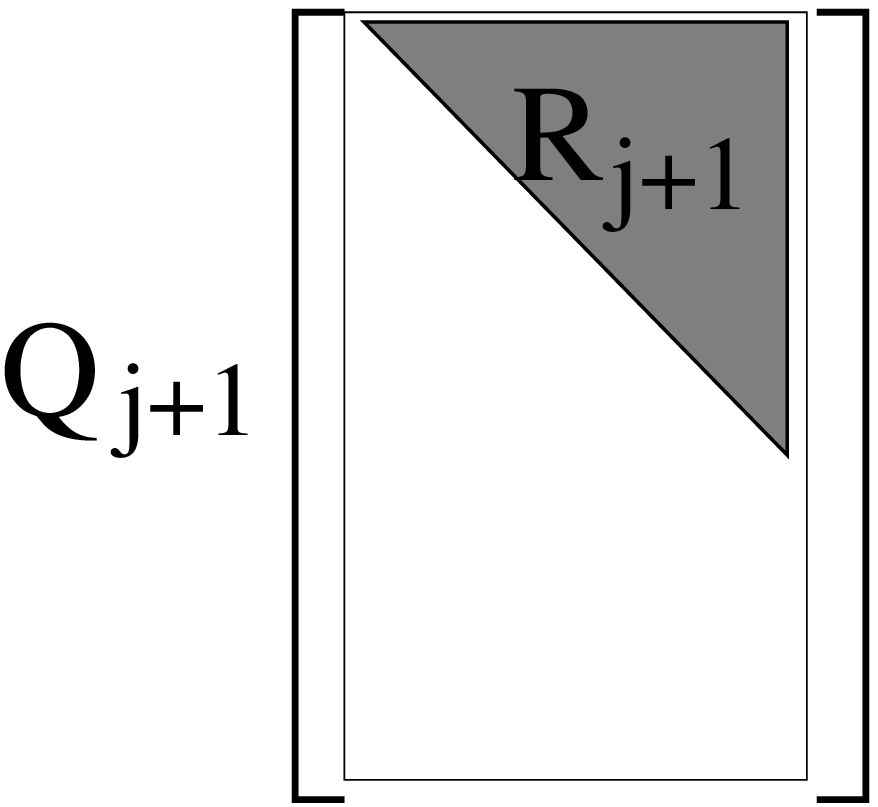,height=1.0in,width=1.3in}}
}}
\label{qrup}
\end{equation}
%
The total flops for decomposing all $S_j(p,q)P_j=Q_jR_j$'s
is $O(n^3)$.

With the $R \,=\,R_j$ at each QR update step, we apply
the iteration (\ref{rrit}) from a randomly generated initial
vector $\bdz_0$.
Theoretically, there is a zero probability that such $\bdz_0$ is
perpendicular to the singular subspace $\spn\{\bdy\}$ and the condition
for convergence may not be satisfied.
In practical floating point computation, however,
the round-off error quickly destroys this orthogonality and
the iteration {\em always} converges.
Moreover, when $\sgmin\big(S_j(p,q)\big)$ is near zero, the convergence
rate (\ref{rritrate}) in Lemma \ref{svcvrt} is quite fast.
If $\sgmin\big(S_j(p,q)\big)$ is, say,  less than $10^{-3}$ of the
second smallest singular value every iteration step in (\ref{rrit}) will
produce 6 correct digits.
It rarely takes more than 3 to 5 iterations to reach
a near zero $\sgmin\big(S_j(p,q)\big)$.

According to (\ref{sgmin25}) in Lemma \ref{lem:sgmin}, when
$\sgmin\big(S_j(p,q)\big) \,<\,\eps\, \sqrt{m-j+1}$ is reached along with 
the singular vector $\bdy$, the entries of $\bdy$ form $\bdv_0$
and $\bdw_0$ that approximate the coefficients of cofactors
$v$ and $w$ respectively with error bound (\ref{cferr}).
By Lemma \ref{lem:uerr}, an approximation $u_0$ to the numerical GCD
can be obtained from solving the linear system
\begin{equation} \label{usys}
\mbox{\scriptsize $
\left[ \begin{array}{c} C_{j}(v_0) \\ C_{j}(w_0) \end{array}
\right]$} \bdu_0 =
\mbox{\scriptsize $
\left[ \begin{array}{c} \bdp \\ \bdq \end{array} \right]$}
\end{equation}
with error bound (\ref{uhaterr}).

The iteration (\ref{rrit}) is applied at $j=n,n-1,\ldots, k$, if
the numerical GCD degree is $k$.
Each step in (\ref{rrit}) requires $O((n-j)^2)$ flops.
The total cost in calculating the numerical GCD degree
is no higher than $O(n^3)$.
The system (\ref{usys}) costs $O(n^3)$ to solve and it is to be
solved only when a possible numerical GCD is detected.

Notice that, inequality \( \sgmin\big(S_{j}(p,q)\big) \,\leq\,
\eps \sqrt{m-j+1} \)
is not a sufficient condition for the given polynomial pair $(p,q)$
to have a numerical GCD of degree $k=j$ within tolerance $\eps$.
Satisfying this inequality alone does not guarantee the existence of
numerical GCD within $\eps$.
Even if the numerical GCD degree is found by this inequality, the
numerical GCD triplet $(\hat{u},\hat{v},\hat{w})$ may not be accurate
enough.
An iterative refinement below verifies the numerical GCD degree and
refines the numerical GCD and cofactors.

\subsection{Iterative refinement} \label{itref}

After obtaining a possible numerical GCD degree $k$, a GCD manifold
$\cPkmn$ of codimension $k$ is tentatively targeted for seeking
minimum distance
to the given polynomial pair $(p,q)$.
Using the degree $k$, we set the numerical GCD system
$\bdf_\bdh(u,v,w)=[\bt;\,\bdp;\,\bdq]$
where $\bdf_\bdh$ is defined in (\ref{fzb}) with an objective in
finding the least squares solution $(u,v,w)$ with $\dg(u) = k$,
$(uv,uw) \in \cPkmn$, and $\big\|(uv,uw)-(p,q)\big\| = \theta_k(p,q)$.

There are several choices for the scaling vector $\bdh$.
If the numerical GCD is required to be monic, then $\bdh = [0;0;\ldots;1]$.
There is a drawback in this choice: when leading coefficients of
$p$ and $q$ are small, forcing $u$ to be monic may cause its remaining
coefficients to be large in magnitude and creating unbalanced system
(\ref{fzb}).
A random vector as $\bdh$ can be a good choice, although there is a zero
probability that such $\bdh$ would make (\ref{fzb}) singular, or a small
probability that the system is ill-conditioned.
The most preferred choice appears to be a scalar multiple of the initial
approximation $\bdu_0$ of the numerical GCD determined by (\ref{usys}).
This is because $u_0$ is close to the numerical GCD $u$ as ensured by
Lemma~\ref{lem:uerr} when $(p,q)$ ~is near manifold $\cPkmn$, and
coefficient vector $\bdu_0$ cannot be perpendicular to $\bdu$.

With the choice of scaling vector 
$\bdh = \beta \bdu_0$ in (\ref{fzb}) and the initial
approximation $(u_0,v_0,w_0)$ described in \S \ref{sec:ini},
the Gauss-Newton iteration (\ref{gnit3}) is applied and iteration stops
when the distance $\dl_j \equiv \|\,\bdf_\bdh(u_j,v_j,w_j)-
[\bt; \,\bdp;\,\bdq]\,\|$ stops decreasing.

This refinement stage outputs the nearness $\rho = \dl_j$ and
the refined numerical GCD triplet $(u,v,w)$.
If $\rho < \eps$, then the $(u,v,w)$ is certified as numerical
GCD triplet for $p$ and $q$.
On the other hand, if the distance $\rho \ge \eps$, then the numerical
GCD degree $k$ that is tentatively determined in \S \ref{sec:ini} is
incorrect and needs to be adjusted downward by one.

At each Gauss-Newton iteration step in (\ref{gnit3}), it is neither 
desirable nor 
necessary to construct the Moore-Penrose inverse 
$J_\bdh(u_j,v_j,w_j)^+$ in an explicit form.
The new iterate $(u_{j+1},v_{j+1},w_{j+1})$ is obtained via solving a
linear least squares problem
\[ \left\{ \begin{array}{l}
\mbox{Solve \ \ $J_\bdh (u_j,v_j,w_j) (\Dl \bdz) ~=~
\bdf_\bdh(u_j,v_j,w_j) - [\bt; \,\bdp;\,\bdq]$
\ ~for~ \ $\Dl \bdz$} \\
\mbox{Set \ \ \ \ \ \ $[\bdu_{j+1};\bdv_{j+1};\bdw_{j+1}] \;=\;
[\bdu_j;\bdv_j;\bdw_j] - \Dl \bdz$.}
\end{array} \right.
\]
While solving $J_\bdh(u_j,v_j,w_j) (\Dl \bdz) = \bdf_\bdh(u_j,v_j,w_j)
- [\bt; \,\bdp;\,\bdq] $
for its least squares solution, the QR decomposition of
$J_\bdh(u_j,v_j,w_j) = Q\,R$ is obtained.
Upon exiting the iteration (\ref{gnit3}), the final upper triangular matrix
$R$ can replace the $R_j$ in (\ref{rrit}) to calculate the smallest
singular value of $J_\bdh(u_j,v_j,w_j)$.
The reciprocal of this singular value is the GCD condition number of the
polynomial pair $(p,q)$ and the computed numerical GCD triplet $(u,v,w)$.
Calculating the condition number requires negligible flops.

The Gauss-Newton iteration here is crucial for its two-fold purpose:

\begin{itemize} 
\item[] \hspace{-8mm} {\em Verifying the numerical GCD.}
By minimizing the residual
\[ \rho ~~=~~ \|\,\bdf_\bdh(u,v,w)-[\bt; \,\bdp;\,\bdq]\,\|,
\]
the Gauss-Newton iteration either
certifies the numerical GCD triplet $(u,v,w)$ from verifying
$\rho < \eps$, or disqualify $k$ as the numerical GCD degree when
$\rho \ge \eps$.
In the latter case the process of computing $\sgmin\big(S_{j}(p,q)\big)$
needs to be continued for decreased $j$ by one.
\item[] \hspace{-8mm} {\em Refining the numerical GCD triplet.}
The Gauss-Newton
iteration filters out the error bounded by (\ref{cferr}) and (\ref{uhaterr}),
obtaining the numerical GCD to the optimal accuracy bounded by
(\ref{cond}).
\end{itemize}

\subsection{The main algorithm and its convergence theorem}
\label{sec:mal}

In summary, the overall algorithm for finding a numerical GCD of a
polynomial pair within a tolerance $\eps$ is described in the following
pseudo-code, which
contains two exit points.

\begin{itemize} 
\item[] \hspace{-6mm} {\bf Algorithm {\sc uvGCD}}
\item {\bf Input}: Pair $(p,q) \in \cPmn$ with $m \ge n$
,backward nearness tolerance $\eps$.
\item Initialize permutation $P_n = I$ and QR decomposition
$S_n(p,q)P_n = Q_n R_n$,
\item {\bf For} $j=n,n-1,\ldots,1$ {\bf do}
\begin{itemize}
\item Apply iteration (\ref{rrit}) on $R =R_j$ and obtain the 
smallest singular value
$\sgmin\big(S_j(p,q)\big) \equiv \sgmin\big(S_j(p,q)P_j\big)$ and
corresponding singular vector $\bdy = P_j[\bdw_0;\,-\bdv_0]$ of
$S_j(p,q)P_j$.
\item {\bf If} $\sgmin\big(S_j(p,q)\big) \,<\, \eps \sqrt{m-j+1}$
then
\begin{itemize} \parskip2mm 
\item Set GCD degree as $k\,=\,j$, extract $v_0$ and $w_0$ from
$\bdy = P_j[\bdw_0;\,-\bdv_0]$,  and compute an initial approximation
$u_0$ to the numerical GCD by solving $(\ref{usys})$.
\item Set up
$\bdf_\bdh(\cdot,\cdot,\cdot)$, $J_\bdh(\cdot,\cdot,\cdot)$
as in (\ref{fzb}) and
(\ref{jac}) with $k=j$ and the scaling vector $\bdh = \bt \bdu_0$
for ~$\bt=1$.
Apply the Gauss-Newton iteration (\ref{gnit3}) with initial iterate
$(u_0,v_0,w_0)$
and terminate the iteration at the triplet $(u,v,w) = (u_l,v_l,w_l)$ 
when the residual
$\dl_l = \|\bdf_\bdh(u_l,v_l,w_l)-[\bt; \,\bdp;\,\bdq]\|$ stops decreasing;
set $\rho = \dl_l$.
\item {\bf If} $\rho < \eps$, then {\bf break} the do-loop,
{\bf end if}
\end{itemize}
\item[] {\bf end if}
\item Update $S_{j-1}(p,q)P_{j-1} = Q_{j-1}R_{j-1}$
as in (\ref{qrup}).
\end{itemize}
\item[] {\bf end do}
\item {\bf Output} GCD triplet $(u,v,w)$ if $\rho <\eps$, or trivial GCD
triplet $(1,p,q)$ if $\rho \ge \eps$.
\end{itemize}

The following is the Numerical GCD Convergence Theorem for the numerical 
GCD algorithm.
The theorem asserts that Algorithm~{\sc uvGCD} converges
to a numerical GCD
and cofactors that can be arbitrarily accurate if the given polynomial
pair is within a sufficiently small perturbation.

\begin{thm}[Numerical GCD Convergence Theorem] 
\label{prop:conv} 
Let $(\hat{p},\hat{q})$ be any\linebreak polynomial pair in $\cPkmn$.
Then for every $\dl >0$, there is an $\eta >0$ such that,
if input items $(p,q)\in \cPmn$ and $\eps > 0$ satisfy
\[
\|(p,q)-(\hat{p},\hat{q})\| ~<~ \eta
~<~ \eps ~<~  \theta_{k+1}(\hat{p},\hat{q})-\eta,
\]
there is a unique numerical GCD $\GCD_{\eps}(p,q) = \GCD(\tilde{p},
\tilde{q})$
with $(\tilde{p},\tilde{q}) \in \cPkmn$ satisfying
$\|(\tilde{p},\tilde{q})-(\hat{p},\hat{q})\| < \dl$.
Moreover, Algorithm {\sc uvGCD} generates a sequence of polynomial
triplets $(u_j,v_j,w_j)$ satisfying 
\[
\lim_{j\rightarrow \infty} u_j ~~=~~ u \in \GCD_\eps(p,q) 
~~~~\mbox{~~~~~~~~and~~~~~~~~}~~~~
\lim_{j\rightarrow \infty} \big\|(u_j v_j, u_j w_j) -
(\tilde{p},\tilde{q})\big\| ~~=~~ 0
\]
\end{thm}

\prf
Assume $m > n$ without loss of generality.
From $\dg\big(\GCD(\hat{p},\hat{q})\big) = k$, we have
$\theta_{k+1}(\hat{p},\hat{q}) > 0$ and we can
choose a $\eta_1$ with $0 < \eta_1 < \theta_{k+1}(\hat{p},\hat{q})/2$.
If $\|(p,q)-(\hat{p},\hat{q})\|\,<\,\eta_1$ and $\eta_1 < \eps <
\theta_{k+1}(\hat{p},\hat{q})-\eta_1$, then
\[ \sgmin\big(S_{k}(p,q)\big) \;<\; \eta_1\,\sqrt{m-k+1} \;<\;
\eps\,\sqrt{m-k+1} \]
by Lemma \ref{lem:sgmin}
and the Gauss-Newton iteration (\ref{gnit3}) will be initiated at certain
$j \ge k$.
For any $j > k$, the distance $\theta_{j}(p,q) \,\ge\,
\theta_{k+1}(p,q) \,\ge\,\theta_{k+1}(\hat{p},\hat{q})-\eta_1 \,>\,
\eps$.
Consequently the Gauss-Newton iteration either diverges or converges
to a point with residual larger than $\eps$.
As a result, Algorithm {\sc uvGCD} will not be terminated at $j > k$.

From $\|(p,q)-(\hat{p},\hat{q})\|\,<\,\eta_1$ and
$\eta_1 < \eps < \theta_{k+1}(\hat{p},\hat{q})-\eta_1$, $\cPkmn$
is the GCD manifold of highest codimension within $\eps$ of $(p,q)$,
namely 
\[ k = \codim\big(\cPkmn\big) =
\max_{0\le j \le n} \big\{\, \codim\big(\cPjmn\big) \big|
\theta_j(p,q) < \eps \,\big\}.
\]
Clearly, $\theta_k(p,q) \le \eta_1 < \eps$ is attainable
at certain $(\tilde{p},\tilde{q}) \in \overline{\cPkmn}$.
Consequently, the unique numerical GCD ~$\GCD_{\eps}(p,q)$ exists
and is identical to the exact GCD of $(\tilde{p},\tilde{q})$.

Let $\hat{u} \in \GCD(\hat{p},\hat{q})$. 
For any fixed $\bdh \in \C^{k+1}$ with $\bdh^\h \hat{\hat{\bdu}} \ne 0$,
let $(\hat{u},\hat{v},\hat{w})$ be the unique solution to the equation
~$\bdf_\bdh(u,v,w) = [1;\hat\bdp;\hat\bdq]$.
By Lemma~\ref{gnlem}, there is a neighborhood $\Dl$ of 
$(\hat{u},\hat{v},\hat{w})$ and a neighborhood $\Sigma$ of 
$(\hat{p},\hat{q})$ such that for all ~$(p,q) \in \Sigma$, the Gauss-Newton
iteration on the system $\bdf_\bdh(u,v,w) = [1,p,q]$ converge to the
least squares solution $(u_*,v_*,w_*)$ from any initial iterate 
~$(u_0,v_0,w_0) \in \Dl$. 
By Lemma~\ref{prop:tricont}, there is an $\eta_2>0$ ~such that 
$(u_0,v_0,w_0) \in \Dl$ and $(p,q) \in \Sigma$
whenever $\|(p,q)-(\hat{p},\hat{q})\| < \eta_2$. 
Set $\eta = \min\{\eta_1,\eta_2\}$, the conclusion of the theorem follows.
\qed

\section{Computing experiment and benchmark}
\label{sec:res}

Our method is implemented as a package {\sc uvGCD}
in Maple and Matlab.
~In addition to a symbolic GCD-finder {\sc gcd},
there are three numerical GCD finders in the SNAP package \cite{jean-lab}
in Maple: {\sc QuasiGCD} \cite{beck-fast}, {\sc EpsilonGCD} \cite{beck-fast},
and {\sc QRGCD} \cite{cor-watt-zhi}.
~Among them {\sc QRGCD} is clearly superior to the other two by a wide margin.
~We thereby compare {\sc uvGCD} with {\sc QRGCD} and {\sc gcd} only.
~Actually, {\sc QuasiGCD} and {\sc EpsilonGCD} output failure messages for all
the test examples in this section.

All test results are obtained on a desktop PC with an Intel Pentium 4 CPU
of 1.8 MHz and 512 Mb memory.
~Unless mentioned specifically (Example 4), both {\sc uvGCD} and {\sc QRGCD}
are tested in Maple 9 with precision set
to 16 digits to simulate hardware precision.

We believe that numerical GCD finders should be tested and compared
based on results from the following aspects.

\begin{enumerate} \parskip-0.5mm
\item Performance on polynomials with increasing numerical GCD
sensitivity.
\item Performance on polynomial having different numerical GCD's
within different tolerance.
\item Performance on numerical GCD's of large degrees.
\item Performance on polynomials with large variation
in coefficient magnitudes.
\item Performance in finding the numerical GCD of ~$(p,p')$~
when ~$p$~ has roots of high multiplicities.
\end{enumerate}

We have established a test suite that includes polynomials
satisfying the above requirements along with those
collected from the literature.
~We demonstrate the robustness and accuracy of {\sc uvGCD} with sample
results below.

{\bf Test 1: A high sensitivity case.} For an even number ~$n$~ and
~$k=n/2$, ~let ~$p_n = u_n\,v_n$~ and ~$q_n = u_n\,w_n$, ~where
\small
\[
u_n = \prod_{j=1}^{k} \left[ \left(x-r_1\al_j\right)^2 +
r_1^2 \bt_j^2 \right], \;
v_n = \prod_{j=1}^{k} \left[ \left(x-r_2\al_j\right)^2 +
r_2^2 \bt_j^2 \right], \] 
\[
w_n = \prod_{j=k+1}^{n} \left[ \left(x-r_1\al_j\right)^2 +
r_1^2 \bt_j^2 \right],
\;\;\;
\al_j =  \cos \frac{j \pi}{n},
\bt_j =  \sin \frac{j \pi}{n}
\] \normalsize
for ~$r_1=0.5$, ~$r_2=1.5$.
~The roots of ~$p_n$~ and ~$q_n$~ spread on the circles of radius ~$0.5$~
and ~$1.5$.
~When ~$n$~ increases, the GCD condition number grows quickly.
~Table \ref{ex1} shows that error on the computed numerical GCD.

\begin{table}[htb] \small
\begin{center}
\begin{tabular}{|l|c|c|c|}
\hline
\ \ $n$ \ \ & condition & QRGCD & uvGCD\\
& number & error & error\\ \hline
$n=6$ & 566.13 & $0.55\times 10^{-14}$ & $0.15\times 10^{-14}$ \\
$n=10$ & 742560.0 & $0.18\times 10^{-11}$ & $0.47\times 10^{-12}$ \\
$n=16$ & $0.33 \times 10^{11}$ & $0.18\times 10^{-4}$ & $0.65\times 10^{-9}$
\\
$n=18$ & $0.17 \times 10^{13}$ & FAIL & $0.53\times 10^{-5}$ \\
$n=20$ & $0.71 \times 10^{14}$ & FAIL & $0.99\times 10^{-6}$ \\ \hline
\end{tabular}
\end{center} 
\caption{Comparison in Test 1} \label{ex1}
\end{table}

{\bf Test 2: Multiple numerical GCD's}. Let
\[ \begin{array}{rcl}
p(x) & = & \prod_{j=1}^{10} ( x - x_j ),  \mbox{\ \ with \ \ }
x_j = (-1)^j \left(\frac{j}{2}\right)
\\
q(x) & = & \prod_{j=1}^{10} \left[ x - x_j + 10^{-j}
\right]  \end{array}
\]
The roots of ~$q$~ have decreasing distances ~$0.1,0.01,\ldots$~ with those
of ~$p$.
~Therefore there are different numerical GCD's  for different tolerances.
~As shown in Table \ref{exam2},
{\sc uvGCD} accurately separates the numerical GCD factors according to
the given tolerance on the listed cases.

\begin{table}[htb]
\small
\begin{center}
\begin{tabular}{|c|c|c|} \hline
tolerance & \multicolumn{2}{c|}{degree (\& nearness) of numerical
GCD found by} \\
\cline{2-3}
$\eps$ & {\sc QRGCD} & {\sc uvGCD}
\\ \hline
$10^{-2}$ & 7 (0.0174) & 9 \ (0.56E-02)
\\
$10^{-3}$ & Fail & 8 \ (0.26E-03)
\\
$10^{-4}$ & Fail & 7 \ (0.14E-04)
\\
$10^{-5}$ & Fail & 6 \ (0.11E-05)
\\
$10^{-6}$ & Fail & 5 \ (0.41E-07)
\\
$10^{-8}$ & Fail & 4 \ (0.42E-08)
\\
$10^{-9}$ & Fail & 3 \ (0.14E-09)
\\
$10^{-10}$ & Fail & 2 \ (0.24E-10)
\\ \hline
\end{tabular}
\end{center} 
\caption{\scriptsize The calculated degrees (and nearness in parentheses)
of numerical GCD within various tolerance on Test 2.} \label{ex2}
\end{table}

{\bf Test 3: numerical GCD of large degrees}. ~For fixed cofactors
~$v(x) = \sum_{j=0}^3 x^j$ ~and ~$w(x) = \sum_{j=0}^4 (-x)^j$, ~let
~$p_n = u_n \,v$~ and ~$q_n = u_n \, w$~ with ~$u_n$~ being a polynomial
of degree ~$n$~ of random integer coefficients in ~$[-5,5]$.
~For the sequence of polynomial pairs ~$(p_n,q_n)$, ~the GCD is
known to be ~$u_n$~ and we can calculate the actual accuracy.
~As shown in Table \ref{ex3}, {\sc uvGCD} maintains its robustness
and high accuracy even when for ~$n$~ reaches ~$2000$, ~while
{\sc QRGCD} works for ~$n < 100$.

\begin{table}[htb] \small
\begin{center}
\begin{tabular}{|c|c|c|} \hline
GCD  & \multicolumn{2}{c|}{coefficient-wise error on computed numerical GCD}
\\ \cline{2-3}
degree & \ \ \ \ \ {\sc QRGCD} \ \ \ \ \  & {\sc uvGCD} \\ \hline
$n=\;\; 50$ & 0.168E-12 & 0.500E-15 ~$\;\,$~ \\
$n=\;\; 80$ & 0.927E-12 & 0.805E-15 ~$\;\,$~ \\
  $n=\;100$ & Fail      & 0.341E-15 ~$\;\,$~ \\
  $n=\;200$ & Fail      & 0.100E-14 ~$\;\,$~ \\
  $n=\;500$ & Fail (*)      & 0.133E-14 ~$\;\,$~ \\
  $n=1000$ & Fail  (*)    & 0.178E-14 ~$\;\,$~ \\
  $n=2000$ & Fail  (*)    & 0.178E-14 ~$\;\,$~ 
\\ \hline
\end{tabular}
\end{center} 
\label{ex3}
\caption{\scriptsize Coefficient errors on random numerical GCD's of degree $n$.
(*): Presumably failed after running hours without results.
}
\end{table}

{\bf Test 4: A case where computing numerical GCD by {\sc uvGCD} is
faster than calculating GCD by Maple.}
~For polynomials with integer coefficients, Maple's symbolic
GCD finder is often faster than {\sc uvGCD}.
~However, {\sc uvGCD} can be substantially more efficient in other cases.
~Here is an example.
~For fixed cofactors ~$v$~ and ~$w$~ as in Test 3, let ~$u_n$~ be
the polynomial of degree ~$n$~ with random rational coefficients and
\begin{equation} \label{fastpair} p_n = u_n \, v,, \;\;\;\; q_n = u_n \, w.
\end{equation}
The GCD is a multiple of ~$u_n$.
~We compare the Maple {\sc gcd} on exact coefficients with our Matlab
{\sc uvGCD} on approximate coefficients.
~Table  \ref{fast} shows the running time on increasing ~$n$.
~In this polynomial series, not only {\sc uvGCD} is faster, the
speed ratio of {\sc uvGCD} over {\sc gcd} increases from ~$2$~
to ~$11$~ when ~$n$~ increases from ~$50$~ to ~$2000$.
~Of course, this result should be taken with caution because
Maple {\sc GCD} always has zero error.

\begin{table}[htb] \small
\begin{center}
\begin{tabular}{||r||r|c||r|c||} \hline
& \multicolumn{2}{|c||}{Maple {\sc gcd}} &
\multicolumn{2}{c||}{\sc uvGCD} \\ \cline{2-5}
 & time & error & time & error \\ \hline
$n = 50$ & 0.25 & 0 & 0.125 & 3.53e-15 \\
$n = 200$ &  7.47 & 0 & 2.437 & 8.69e-14 \\
$n = 1000$ & 574.90 & 0 & 82.270 & 1.64e-13 \\
$n = 2000$ & 10910.60 & 0 & 969.625 & 1.82e-12 \\
\hline
\end{tabular}
\end{center} 
\caption{\scriptsize Comparison between Maple's symbolic {\sc gcd}
and {\sc uvGCD} on polynomial pairs ~$(p_n,q_n)$~ in (\ref{fastpair})}
\label{fast}
\end{table}

\begin{figure}[t]
\begin{center}
\epsfig{figure=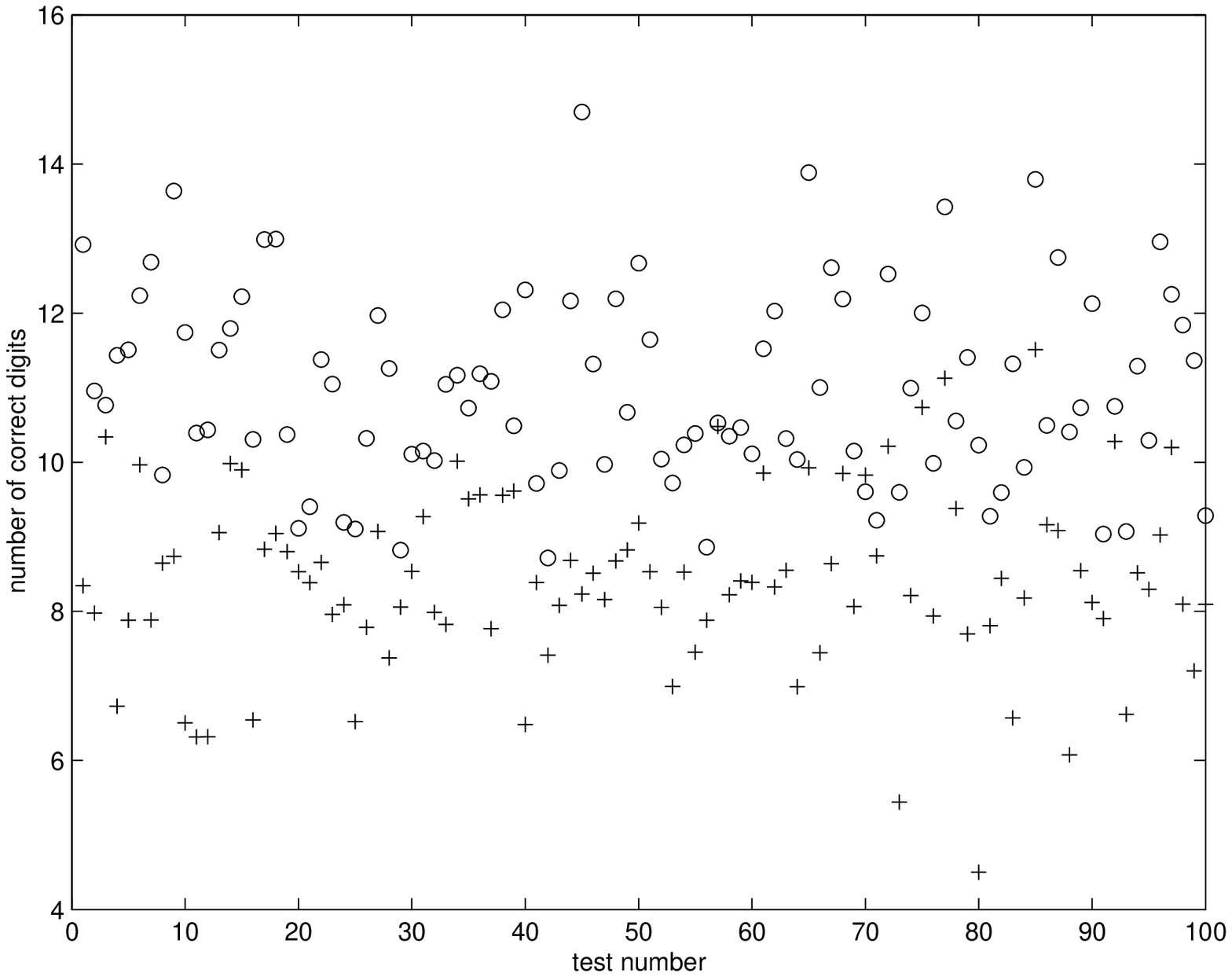,height=2.2in,width=3.6in}
\end{center} 
\centerline{\scriptsize {\tiny $\bigcirc$}: {\sc uvGCD}, \ \ \ \ \ \ \
$+$: {\sc QRGCD} }
\caption{\scriptsize Accuracy comparison on 100 polynomial pairs in Test 4.
~The height of each point is the number of correct digits
approximating the numerical GCD coefficients by {\sc uvGCD} or {\sc QRGCD}
at each polynomial pair.}\label{jmp1}
\end{figure}

{\bf Test 5: Numerical GCD with large variation in coefficient magnitudes}.
~For fixed ~$v$~ and ~$w$~ as in Test 3, let
\[ u(x) = \sum_{j=0}^{15} c_{j} \, 10^{e_{j}} x^j
\]
where for every ~$j$,  ~$c_{j}$~ and ~$e_{j}$~
are random integers in ~$[-5,5]$~ and ~$[0,6]$~ respectively.
~The polynomial pair ~$p = u\,v$~ and ~$q = u\,w$~ are then constructed while
{\sc QRGCD} and {\sc uvGCD} are called to find the numerical GCD of ~$(p,q)$.
Notice that ~$u$~ is the known GCD whose coefficient jumps between ~$0$~ and
~$5\times 10^6$~ in magnitude.
After applying the numerical GCD finders on each pair ~$(p,q)$,
~we calculated the coefficient-wise relative errors ~$\theta$~ and
~$\vartheta$~ of {\sc QRGCD} and {\sc uvGCD} respectively.
Roughly speaking, ~$-\log_{10}\theta$~ and ~$-\log_{10}\vartheta$~ are the
minimum number of correct digits obtained for approximating coefficients of
~$u$~ by {\sc QRGCD} and {\sc uvGCD} respectively.
This test is repeated 100 times.
~Figure \ref{jmp1} shows that on average {\sc QRGCD} gets about 8 digits
correct on each coefficient, while {\sc uvGCD} attains about 11.

\begin{figure}[htb]
\begin{center}
\epsfig{figure=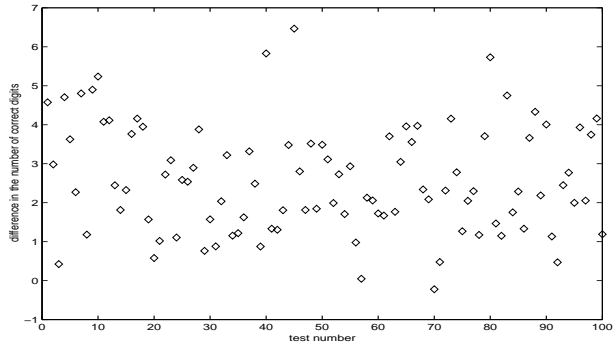,height=1.8in,width=3.2in}
\end{center} 
\caption{\scriptsize Accuracy comparison on 100 polynomial pairs in Test 4.
~The vertical axis is the difference in the number of correct digits
approximating the numerical GCD coefficients by {\sc uvGCD} or {\sc QRGCD}
at each test.}\label{jmp2}
\end{figure}

Figure \ref{jmp2} shows the {\em difference} in the number of correct digits
obtained on coefficients from each test.
~On those 100 tests, {\sc uvGCD} obtains up to 6.5 more correct
digits than {\sc QRGCD} on 99 test, while slightly less accurate than
{\sc QRGCD} on only one polynomial pairs (i.e. the test 70).

{\bf Test 6: GCD of ~$p$~ and ~$p'$.} ~Let ~$p$~ be
\[
p(x) = (x-1)^{m_1} (x-2)^{m_2} (x-3)^{m_3} (x-4)^{m_4} \]
for different sets of ~$m_1,m_2,m_3,m_4$.
~Finding the numerical GCD of ~$p$~ and ~$p'$~ may be difficult for some
numerical GCD finders,
as shown in Table \ref{ex5} for {\sc QRGCD} and \cite{rupprecht}.
~This numerical GCD computation has an important application in polynomial
root-finding.
~On the other hand, {\sc uvGCD} is originally built for this purpose
and shows its tremendous robustness.

\begin{table}[ht]
\small
\begin{center}
\begin{tabular}{|r|c|c|c|} \hline
& \multicolumn{3}{c|}{coefficient-wise relative error} \\ \cline{2-4}
$[m_1,m_2,m_3,m_4]$ & Maple & Maple &   \\
& {\sc QRGCD} & {\sc gcd} & {\sc uvGCD} \\ \hline
$[2,1,1,0]$ & 1.0E-13 & 1.0e-16 & 6.7E-16 \\
$[3,2,1,0]$ & 1.5E-12 & 1.0e-16 & 1.8E-14 \\
$[4,3,2,1]$ & 1.6E-07 & 1.0e-16 & 4.5E-14 \\
$[5,3,2,1]$ & Fail & 3.5e-16 & 4.6E-13 \\
$[9,6,4,2]$ & Fail & Fail(*) & 3.5E-12 \\
$[20,14,10,5]$ & Fail & Fail(*) & 1.7E-12 \\
$[80,60,40,20]$ & Fail & Fail(*) & 3.5E-11 \\
$[100,60,40,20]$ & Fail & Fail(*) & 2.6E-11 \\
\hline
\end{tabular}
\end{center} 
\caption{\scriptsize
Comparison on $\GCD(p,p')$ for $p$
in Test 5.
\newline (*): Symbolic {\sc gcd} fails because $p$ is no longer exact.}
\label{ex5}
\end{table}

\bibliographystyle{amsalpha}

\begin{thebibliography}{A}

\bibitem{barnett}
S.~Barnett.
\newblock {\em {~Polynomials} and Linear Control Systems}.
\newblock ~Monographs and textbooks in pure and applied mathematics, Marcel
  Dekker, Inc, New York, 1983.

\bibitem{beck-fast}
B.~Beckermann and G.~Labahn.
\newblock {~A} fast and numerically stable {E}uclidean-like algorithm for
  detecting relatively prime numerical polynomials.
\newblock {\em ~J. Symb. Comp.}, 26:691--714, 1998.

\bibitem{bre-kun}
R.~P. Brent and H.~T. Kung.
\newblock {~Systolic} {VLSI} arrays for polynomial {GCD} computation.
\newblock {\em ~IEEE Trans. on Computers}, C-33:731--736, 1984.

\bibitem{chin-corless}
P.~Chin, R.~M. Corless, and G.~F. Corless.
\newblock {~Optimization} strategies for the approximate {GCD} problem.
\newblock ~Proc. ISSAC '98, ACM Press, pp 228-235, 1998.

\bibitem{chousc}
S.-C. Chou.
\newblock {\em {~Mechanical} Geometry Theorem Proving}.
\newblock ~D. Reidel Publishing Co., Dordrecht Holand, 1988.

\bibitem{corless-gianni}
R.~M. Corless, P.~M. Gianni, B.~M. Trager, and S.~M. Watt.
\newblock {~The} singular value decomposition for polynomial systems.
\newblock ~Proc. ISSAC '95, ACM Press, pp 195-207, 1995.

\bibitem{cor-watt-zhi}
R.~M. Corless, S.~M. Watt, and L.~Zhi.
\newblock {~{QR}} factoring to compute the {GCD} of univariate approximate
  polynomials.
\newblock {\em ~IEEE Trans. Signal Processing}, 52:3394--3402, 2003.

\bibitem{DemmelBook}
J.~W. Demmel.
\newblock {\em ~{Applied} Numerical Linear Algebra}.
\newblock ~SIAM, Philadelphia, 1997.

\bibitem{dun}
D.~K. Dunaway.
\newblock {~Calculation} of zeros of a real polynomial through factorization
  using {Euclid} algorithm.
\newblock {\em ~SIAM J. Numer. Anal.}, 11:1087--1104, 1974.

\bibitem{emiris-galligo-lombardi}
I.~Z. Emiris, A.~Galligo, and H.~Lombardi.
\newblock {~Certified} approximate univariate {GCD}s.
\newblock {\em ~J. Pure Appl. Algebra}, 117/118:229--251, 1997.

\bibitem{fau93}
O.~D. Faugeras.
\newblock {\em {~Three}-Dimentional Computer Vision: A Geometric Viewpoint}.
\newblock ~MIT Press, Cambridge, Mass., 1993.

\bibitem{GatGer}
J.~V. Gathern and J.~Gerhard.
\newblock {\em {~Modern} Computer Algebra}.
\newblock ~Cambridge University Press, second edition, 2003.

\bibitem{geddes}
K.~O. Geddes, S.~R. Czapor, and G.~Labahn.
\newblock {\em {~Algorithms} for Computer Algebra}.
\newblock ~Kluwer Academic Publishers, Boston, 1992.

\bibitem{golub-vanloan}
G.~H. Golub and C.~F. {Van Loan}.
\newblock {\em {~Matrix} Computations}.
\newblock ~The John Hopkins University Press, Baltimore and London, 3rd
  edition, 1996.

\bibitem{hen-seb}
D.~Henrion and M.~Sebek.
\newblock {~Reliable} numerical methods for polynomial matrix triangulation.
\newblock {\em ~IEEE Trans. on Automatic Control}, 44:497--501, 1997.

\bibitem{hribernig-stetter}
V.~Hribernig and H.~J. Stetter.
\newblock {~Detection} and validation of clusters of polynomial zeros.
\newblock {\em ~J. Symb. Comput.}, 24:667--681, 1997.

\bibitem{jean-lab}
C.-P. Jeannerod and G.~Labahn.
\newblock {~The} {SNAP} package for arithmetic with numeric polynomials.
\newblock ~In International Congress of Mathematical Software, World
  Scientific, pages 61-71, 2002.

\bibitem{kai-nod}
H.~Kai and M.-T. Noda.
\newblock {~Hybrid} rational approximation and its applications.
\newblock {\em ~Reliable Computing}, 6:429--438, 2000.

\bibitem{kmyz05}
E.~Kaltofen, J.~May, Z.~Yang, and L.~Zhi.
\newblock {~Structured} low rank approximation of {Sylvester} matrix.
\newblock {\em ~in Symbolic-Numeric Computation, Trends in Mathematics, D. Wang
  and L. Zhi, editors, Birkh\"auser Verlag, Basel, Switzerland}, pages 69--83,
  2007.

\bibitem{KalYanZhi06}
E.~Kaltofen, Z.~Yang, and L.~Zhi.
\newblock {~Approximate} greatest common divisor of several polynomials with
  linearly constrained coefficients and singular polynomials.
\newblock ~Proc. ISSAC'06, ACM Press, pp 169--176, 2006.

\bibitem{kar-mit-00}
N.~Karcanias and M.~Mitrouli.
\newblock {~Numerical} computation of the {Least Common Multiple} of a set of
  polynomials.
\newblock {\em ~Reliable Computing}, 6:439--457, 2000.

\bibitem{kar-Lak96}
N.~K. Karmarkar and Y.~N. Lakshman.
\newblock {~Approximate} polynomial greatest common divisors and nearest
  singular polynomials.
\newblock ~Proc. ISSAC'96, pp 35-42, ACM Press, 1996.

\bibitem{karmarkar-lakshman}
N.~K. Karmarkar and Y.~N. Lakshman.
\newblock {~On} approximate polynomial greatest common divisors.
\newblock {\em ~J. Symb. Comput.}, 26:653--666, 1998.

\bibitem{li-zeng-03}
T.-Y. Li and Z.~Zeng.
\newblock ~{A} rank-revealing method with updating, downdating and
  applications.
\newblock {\em SIAM J. Matrix Anal. Appl.}, 26:918--946, 2005.

\bibitem{bjorck96}
\mbox{\AA}ke Bj\"{o}rck.
\newblock {\em {~Numerical} Methods for Least Squares Problems}.
\newblock ~SIAM, Philadelphia, 1996.

\bibitem{mer90}
J.~P. Merlet.
\newblock {\em {~Les} Robots Parall\`eles}.
\newblock ~Trait\'es de Nouvellles Technologiques, Herme\`es, 1990.

\bibitem{noda-sas}
M.-T. Noda and T.~Sasaki.
\newblock {~Approximate} {GCD} and its application to ill-conditioned algebraic
  equations.
\newblock {\em ~J. Comput. Appl. Math.}, 38:335--351, 1991.

\bibitem{pan96}
V.~Y. Pan.
\newblock {~Numerical} computation of a polynomial {GCD} and extensions.
\newblock {\em ~Information and Computation}, 167:71--85, 2001.

\bibitem{lia-pil}
S.~Pillai and B.~Liang.
\newblock {~Blind} image deconvolution using {GCD} approach.
\newblock {\em ~IEEE Trans. Image Processing}, 8:202--219, 1999.

\bibitem{rupprecht}
D.~Rupprecht.
\newblock {~An} algorithm for computing certified approximate {GCD} of n
  univariate polynomials.
\newblock {\em ~J. Pure and Appl. Alg.}, 139:255--284, 1999.

\bibitem{schonhage}
A.~Sch\"onhage.
\newblock {~Quasi}-{GCD} computations.
\newblock {\em ~J. Complexity}, 1:118--137, 1985.

\bibitem{sed-chang}
T.~W. Sederberg and G.~Z. Chang.
\newblock {~Best} linear common divisors for approximate degree reduction.
\newblock {\em ~Comp.-Aided Design}, 25:163--168, 1993.

\bibitem{stetter_book}
H.~J. Stetter.
\newblock {\em {~Numerical} Polynomial Algebra}.
\newblock ~SIAM, 2004.

\bibitem{sto-so}
P.~Stoica and T.~S\"oderstr\"om.
\newblock {~Common} factor detection and estimation.
\newblock {\em ~Automatica}, 33:985--989, 1997.

\bibitem{Tay00}
J.~Taylor.
\newblock {\em {~Several} Complex Variables with Connections to Algebraic
  Geometry and Lie Groups}.
\newblock ~American Mathematical Society, Providence, Rhode Island, 2000.

\bibitem{ver-wang}
J.~Verschelde and Y.~Wang.
\newblock {~Computing} dynamic output feedback laws.
\newblock {\em ~IEEE Trans. Automatic Control}, pages 1552--1571, 2004.

\bibitem{WinAll08b}
J.~R. Winkler and J.~D. Allan.
\newblock {~Structured} low rank approximations of the {Sylvester} resultant
  matrix for approximate {GCDs} of {Bernstein} basis polynomial.
\newblock {\em ~Electronic Transactions on Numerical Analysis}, pages 141--155,
  2008.

\bibitem{WinAll08a}
J.~R. Winkler and J.~D. Allan.
\newblock {~Structured} total least norm and approximate {GCD}s of inexact
  polynomial.
\newblock {\em ~J. of Computational and Applied Math.}, pages 1--13, 2008.

\bibitem{zaro}
C.~J. Zarowski, X.~Ma, and F.~W. Fairman.
\newblock {~A} {QR}-factorization method for computing the greatest common
  divisor of polynomials with real-valued coefficients.
\newblock {\em ~IEEE Trans. Signal Processing}, 48:3042--3051, 2000.

\bibitem{zeng_multroot}
Z.~Zeng.
\newblock {~Algorithm} 835: {MultRoot} -- {A} {Matlab} package for computing
  polynomial roots and multiplicities.
\newblock {\em ACM Trans. Math. Software}, 30:218--235, 2004.

\bibitem{zeng-mr-05}
Z.~Zeng.
\newblock {~Computing} multiple roots of inexact polynomials.
\newblock {\em ~Math. Comp.}, 74:869--903, 2005.

\bibitem{apatools}
Z.~Zeng.
\newblock ~{ApaTools}: A {Maple} and {Matlab} toolbox for approximate
  polynomial algebra.
\newblock In M.~Stillman, N.~Takayama, and J.~Verschelde, editors, {\em
  Software for Algebraic Geometry, IMA Volume 148}, pages 149--167. ~Springer,
  2008.

\bibitem{ZengNpa}
Z.~Zeng.
\newblock {~Regularization} and matrix computation in numerical polynomial
  algebra.
\newblock {~In} {\em Approximate Commutative Algebra}, SpringerWienNewYork, L.
  Robbinano and J. Abbott eds., pp. 125--162, 2009.

\bibitem{ZengAIF}
Z.~Zeng.
\newblock {~The} approximate irreducible factorization of a univariate
  polynomial. ~{Revisited}.
\newblock ~Proceedings of ISSAC '09, ACM Press, pp. 367--374, 2009.

\bibitem{zeng-dayton}
Z.~Zeng and B.~Dayton.
\newblock {~The} approximate {GCD} of inexact polynomials. {II: A} multivariate
  algorithm.
\newblock ~Proceedings of ISSAC'04, ACM Press, pp 320-327, (2004).

\bibitem{zeng-li-jcf}
Z.~Zeng and T.-Y.~Li.
\newblock {~A} numerical method for computing the {Jordan Canonical Form}.
\newblock ~Preprint, 2007.

\bibitem{zippel}
R.~Zippel.
\newblock {\em {~Effective} Polynomial Computation}.
\newblock ~Kluwer Academic Publishers, Boston, 1993.

\end{thebibliography}

\end{document}